\documentclass[draft]{amsart}
\usepackage{times}
\usepackage{amsfonts}
\usepackage{amsmath}
\usepackage{amscd}
\usepackage{amssymb}
\usepackage{latexsym}
\usepackage{amsthm}

\theoremstyle{plain}
\newtheorem{stelling}{Theorem}

\swapnumbers
\newtheorem{theorem}[subsection]{Theorem}
\newtheorem{corollary}[subsection]{Corollary}
\newtheorem{lemma}[subsection]{Lemma}
\newtheorem{proposition}[subsection]{Proposition}

\theoremstyle{definition}
\newtheorem{definition}[subsection]{Definition}

\theoremstyle{remark}

\newtheorem{remark}[subsection]{Remark}

\swapnumbers

\newcommand\preprintnote {preprint on \myhomepage}
\newcommand\myhomepage{http://www.math.ohio-state.edu/\-\~{}schoutens}

\newcommand{\emptyprop}{q}

\newcommand \id{\mathfrak a}
\newcommand \into{\hookrightarrow}
\newcommand \inv[1]{{#1^{-1}}}
\newcommand \inverse[2]{{#1^{-1}(#2)}}
\newcommand \iso{\cong}
\newcommand \loc{{\mathcal {O}}}
\newcommand \map[1]{{\newcommand{\tmpprop}{#1q}  \if\tmpprop\emptyprop \to\else \xrightarrow{{\phantom{i}{#1}\phantom{i}}}\fi}} 
\newcommand \maxim{\mathfrak m}
\newcommand \nat{\mathbb N}

\newcommand \pol[2]{#1[#2]}

\newcommand \pr{\mathfrak p}

\newcommand \restrict [2]{\left.#1\right|_{{#2}}}
\newcommand \rij[2]{(#1_1,\dots,#1_{#2})}

\newcommand \sheaf{{\mathcal {F}}}
\newcommand \tensor{\otimes}
\newcommand \tor[4]{\operatorname{Tor}^{#1}_{#2}(#3,#4)}
\newcommand \op\operatorname

\newcommand{\commdiagram}[9][]{%
\begin{equation}
{\newcommand{\tmpprop}{#1q} 
\if\tmpprop\emptyprop \relax\else \label{#1}\fi}
\begin{aligned}%
\mbox{
\begin{picture}(130,90)%
\put(120,70){\vector( 0,-1){50}}%
\put(10,80){\vector( 1, 0){100}}%
\put(0,70){\vector( 0,-1){50}}%
\put(10,10){\vector( 1, 0){100}}%
\put(115,80){\makebox(0,0)[l]{$#4$}}%
\put(5,80){\makebox(0,0)[r]{$#2$}}%
\put(115,10){\makebox(0,0)[l]{$#9$}}%
\put(5,10){\makebox(0,0)[r]{$#7$}}%
\put(-3,50){\makebox(0,0)[r]{$#5$}}
\put(123,50){\makebox(0,0)[l]{$#6$}}
\put(60,3){\makebox(0,0)[c]{$#8$}}
\put(60,88){\makebox(0,0)[c]{$#3$}}
\end{picture}}
\end{aligned}
\end{equation}}

\newcommand\commtrianglefront[7][]{%
\begin{equation}
{\newcommand{\tmpprop}{#1q} 
\if\tmpprop\emptyprop \relax\else \label{#1}\fi}
\begin{aligned}%
\mbox{
\begin{picture}(120,80)%
\put(55,70){\vector(-1,-2){30}}
\put(65,70){\vector(1,-2){30}}
\put(30,5){\vector(1,0){60}}
\put(60,75){\makebox(0,0)[c]{$#2$}}
\put(25,5){\makebox(0,0)[r]{$#4$}}
\put(95,5){\makebox(0,0)[l]{$#6$}}
\put(60,0){\makebox(0,0)[c]{$#5$}}
\put(37,43){\makebox(0,0)[r]{$#3$}}
\put(83,43){\makebox(0,0)[l]{$#7$}}
\end{picture}}
\end{aligned}
\end{equation}}

\newcommand\commtriangleback[7][]{%
\begin{equation}
{\newcommand{\tmpprop}{#1q}
\if\tmpprop\emptyprop \relax\else \label{#1}\fi}
\begin{aligned}%
\mbox{
\begin{picture}(120,80)%
\put(55,70){\vector(-1,-2){30}}
\put(65,70){\vector(1,-2){30}}
\put(30,5){\vector(1,0){60}}
\put(60,75){\makebox(0,0)[c]{$#2$}}
\put(25,5){\makebox(0,0)[r]{$#6$}}
\put(95,5){\makebox(0,0)[l]{$#4$}}
\put(60,0){\makebox(0,0)[c]{$#5$}}
\put(37,43){\makebox(0,0)[r]{$#7$}}
\put(83,43){\makebox(0,0)[l]{$#3$}}
\end{picture}}
\end{aligned}
\end{equation}}

\newcommand \acf{algebraically closed field}

\newcommand \ch{characteristic}
\newcommand \homo{homomorphism}
\newcommand \CM{Coh\-en-Mac\-au\-lay}
\renewcommand\iff{if, and only if,}

\DeclareMathOperator*{\UP}{ulim}
\newcommand \up[1]{\UP_{#1\to\infty}}
\newcommand \ul[1]{\seq{#1}\infty}
\newcommand \uf{\varphi}
\newcommand \seq[2]{#1\mathstrut_{#2}}
\newcommand \sect[2]{\langle#2\rangle_{#1}}
\newcommand \multf[3]{{#1#2}_{#3}}
\newcommand \sr{approximation}
\newcommand \Sr{Approximation}
\newcommand \BS{Brian\c{c}on-Skoda}

\newcommand  \nsreg{weakly  difference regular}
\newcommand  \nsregglob{difference regular}
\newcommand  \genrat{generically F-rational}
\newcommand  \genreg{weakly generically F-regular}
\newcommand  \genregglob{generically F-regular}
\newcommand  \qreg{ultra-F-regular}
\newcommand  \qpure{ultra-F-pure}
\newcommand  \qpurity{ultra-F-purity}
\newcommand  \qfrob{ultra-Frobenius}

\newcommand  \Qfrobs{Ultra-Frobenii}
\newcommand  \Cech{\u{C}ech}
\newcommand  \los{\L os' Theorem}
\newcommand  \BCM[1]{\mathcal B(#1)}

\newcommand \name[1]{#1}

\newcommand  \bc[1]{{#1}_*}
\newcommand \frob[1]{\uf_{#1}}
\newcommand \ulfrob{\frob\infty}

\newcommand \polar[1]{#1^\lozenge}

\newcommand \grad[2]{\left [#1\right]_{#2}}
\newcommand \gr[1]{#1^{\text{gr}}}

\newcommand \ac[2]{{\mathbb{#1}_{#2}^{\text{alg}}}}
\newcommand \zet{\mathbb Z}
\newcommand \cc{\mathcal C^{\bullet}}

\title {Log-terminal Singularities and Vanishing Theorems}
\author{Hans Schoutens}
\thanks{Partially supported by a  grant from the National Science Foundation and by visiting positions at Paris VII and at the Ecole Normale Sup\'erieure.}
\date{15.03.2003}
\address{Department of Mathematics\\
Ohio State University\\
Columbus, OH 43210 (USA)}
\email{schoutens@math.ohio-state.edu}
\urladdr{\myhomepage} 

\begin{document}

\begin{abstract}   
Generalizing work of Smith and Hara, we give a new characterization of log-terminal singularities for finitely generated algebras over $\mathbb C$, in terms of purity properties of  ultraproducts of characteristic $p$ Frobenii. 

The first application is  a Bout\^ot-type theorem  for log-terminal singularities: given a pure morphism  $Y\to X$ between affine $\mathbb Q$-Gorenstein varieties of finite type over $\mathbb C$, if $Y$ has at most a log-terminal singularities, then so does $X$. The second application is the Vanishing for Maps of Tor for log-terminal singularities: if $A\subset R$ is a Noether Normalization of a finitely generated $\mathbb C$-algebra $R$ and  $S$ is a  finitely generated $R$-algebra with log-terminal singularities, then the natural morphism  $\operatorname{Tor}^A_i(M,R) \to \operatorname{Tor}^A_i(M,S)$ is zero, for every $A$-module $M$ and every $i\geq 1$. The final application is the Kawamata-Viehweg Vanishing Theorem for  a connected projective variety $X$ of finite type over $\mathbb C$ whose affine cone has a log-terminal vertex (for some choice of polarization). As a smooth Fano variety has this latter property, we obtain a proof of the following conjecture  of  Smith for quotients of smooth Fano varieties:  if $G$ is the complexification of a real Lie group acting algebraically on a projective smooth Fano variety $X$, then for any numerically effective line bundle $\mathcal L$ on any GIT quotient $Y:=X/\!/G$, each cohomology module $H^i(Y,\mathcal L)$ vanishes for $i>0$, and, if $\mathcal L$ is moreover big, then $H^i(Y,\mathcal L^{-1})$ vanishes for $i<\operatorname{dim}Y$. 
\end{abstract}

\maketitle

\section{Introduction}

The work of \name{Smith}, \name{Hara} et al., has led to a characterization of log-terminal singularities (equivalence (\ref{i:lt})$\,\Longleftrightarrow\,$(\ref{i:lt}') below) in terms of purity properties of Frobenius  on a general reduction modulo $p$ (\emph{F-regular type}). Although this characterization has proven to be very useful, one of its main drawbacks is the fact that it is not known to be inherited by  quotients of group actions. The main technical result of the present paper is a similar characterization (\emph{\qreg{ity}}) without this defect. 

\begin{stelling}\label{T:MT}
Let $R$ be a local $\mathbb Q$-Gorenstein domain essentially of finite type over a field of \ch\ zero. Then the following are equivalent:
\begin{enumerate}
\item\label{i:lt} $R$ has  log-terminal singularities.
\item[(\ref{i:lt}')] $R$ is F-regular type.
\item\label{i:uFreg} $R$ is \qreg.
\end{enumerate}
\end{stelling}

The implication (\ref{i:lt}')$\implies$(\ref{i:lt}) is proven in \cite[Corollary 4.16]{SmVan} or \cite{HW}, using \name{Smith}'s work on rational singularities in \cite{SmFrat}; the converse implication is proven by \name{Hara} in \cite[Theorem 5.2]{HaRat}. We will give a proof  in \S\ref{s:pfMT} below for the equivalence of \eqref{i:uFreg} with the two other conditions. The notion of `\qreg{ity}' should be viewed as a non-standard version of the notion of `strong F-regularity'. More precisely, let $R$ be a local domain essentially of finite type over $\mathbb C$ (see Remark~\ref{R:bf} below for arbitrary fields). In \cite{SchNSTC}, we associated to $R$ a canonically defined extension $\ul R$, called the \emph{non-standard hull} of $R$, which is realized as the ultraproduct of certain  local rings in \ch\ $p$, called \emph{\sr{s}} of $R$ (see \S\ref{s:affalg} below for exact definitions). One should view  an \sr\ of $R$  as a more canonical way of  reducing $R$ modulo $p$ (see \S\ref{R:sr}), and a non-standard hull of $R$, as a convenient way of storing all these reductions into a single algebraic object. With an \emph\qfrob\  on $R$, we mean the ring \homo\ into the non-standard hull $\ul R$ given by the rule $x\mapsto x^\pi$, where $\pi$ is a non-standard integer obtained as the ultraproduct of various powers of prime numbers (see \S\ref{s:uf} for precise definitions). We call $R$ \emph\qreg, if for each non-zero $c$ in $R$, we can find an \qfrob\ $x\mapsto x^\pi$ such that the $R$-module morphism $R\to \ul R\colon x\mapsto cx^\pi$ is pure. One should compare this with the \name{Hochster-Huneke} notion of  \emph{strong F-regularity} of a domain $R$ of prime \ch\ $p$:   for each non-zero $c$ in $R$, there is  a power $q$ of $p$, such that the morphism $R\to R\colon x\mapsto cx^q$ is split (which under these conditions is equivalent with it being pure). If $R$ is moreover $\mathbb Q$-Gorenstein then strongly F-regular is equivalent by \cite{MacC} with  \emph{weakly F-regular}, that is to say, with the property that every ideal is tightly closed.

\subsection*{Application 1: Quotients of Log-terminal Singularities}
The first application shows that log-terminal singularities are preserved under quotients of reductive groups, provided the quotient is $\mathbb Q$-Gorenstein. Although this seems to be a result that ought to have a proof using Kodaira Vanishing (as for instance in \cite{Bou}), I do not know of any argument other than the one provided here (see  Remark~\ref{R:bou}; similar descent properties for rational singularities, that is to say, for the main result in \cite{Bou}, are treated in \cite{SchBCM,SchRatSing}).

\begin{stelling}\label{T:ltpure}
Let $R\to S$ be a local \homo\ of $\mathbb Q$-Gorenstein local domains essentially of finite type over a field of \ch\ zero. If $R\to S$ is cyclically pure and if  $S$ has log-terminal singularities, then so has $R$. 

In particular, let $G$ be a reductive  group acting algebraically on an affine $\mathbb Q$-Gorenstein variety $X$. If $X$ has at most   log-terminal singularities, then so has the quotient space $X/G$, provided it is  $\mathbb Q$-Gorenstein.
\end{stelling}

\subsection*{Application 2: Vanishing of Maps of Tor}
The next result (see Theorem~\ref{T:vmt} below for  the proof) was previously only known for $S$ regular (\cite[Theorem 9.7]{HuTC}), or more generally, for $S$ weakly CM$\mathstrut^n$-regular (\cite[Theorem 4.12]{HHbigCM2}).

\begin{stelling}\label{T:vantor}
Let $R\to S$ be  a \homo\ of $\mathbb C$-affine algebras such that $S$ is a domain with at most log-terminal singularities (or, more generally, a pure subring of such a ring). Let $A$ be a regular subring of $R$ over which $R$ is module finite. Then for every $A$-module $M$ and every $i\geq 1$, the natural morphism $\tor AiMR\to \tor AiMS$ is zero.
\end{stelling}

\subsection*{Application 3: Vanishing Theorems}
Purity of Frobenius was used effectively in \cite{HR} to prove the \CM{ness} of ring of invariants. Exploiting this further, \name{Mehta} and \name{Ramanathan} deduced Vanishing Theorems for Schubert varieties from purity properties of Frobenius in \cite{MR}. The approach in this paper is a non-standard analogue of these ideas, especially those from \cite{SmFano}. Let $X$ be a connected, normal projective variety of \ch\ zero. Recall that $\op{Spec} S$ is called an \emph{affine cone} of $X$, if $S$ is some finitely generated graded algebra such that  $X=\op{Spec} S$ (for each choice of ample invertible sheaf on $X$, one obtains such a graded ring $S$; see \S\ref{s:pol} below). The \emph{vertex} of the affine cone is by definition the closed point on $\op{Spec}S$ determined by the irrelevant maximal ideal of $S$ (generated by all  homogeneous elements of positive degree). We call $X$ \emph{globally \qreg}, if some affine cone of $X$ has an  \qreg\ vertex. In particular, in view of Theorem~\ref{T:MT}, if the vertex of an affine cone is a log-terminal singularity, then $X$ is globally \qreg. Since the anti-canonical cone of a smooth Fano variety (or more generally, a Fano variety with rational singularities) has this property (see Theorem~\ref{T:Fano} below), every smooth Fano variety is globally \qreg, and, more generally, by Theorem~\ref{T:ltpure}, so is any GIT (Geometric Invariant Theory) quotient   of a smooth Fano variety. In Corollaries~\ref{C:nef} and \ref{C:bignef} we will show the following vanishing of cohomology  for globally \qreg\ varieties.

\begin{stelling}\label{T:KVvan}
Let $X$ be a globally \qreg\ projective variety and let $\mathcal L$  be a numerically effective  line bundle on $X$ (this includes the case $\mathcal L=\loc_X$). For each $i>0$, the cohomology module $H^i(X,\mathcal L)$ vanishes. Moreover, if $\mathcal L$ is also  big, then $H^i(X,\inv{\mathcal L})$ vanishes for each $i<\op{dim}X$.
\end{stelling}

The following particular instance of this theorem was originally conjectured by \name{Smith} in \cite{SmFano}.

\begin{stelling}\label{T:Fanoquot}
Let   $G$ be a reductive group acting algebraically on a projective Fano variety $X$ with rational singularities and let $Y:=X/\!/G$ be a  GIT quotient of $X$ (with respect to some linearization of the action of $G$). If $\mathcal L$ is a numerically effective line bundle  on  $Y$, then each cohomology module $H^i(Y,\mathcal L)$ vanishes for $i>0$, and, if $\mathcal L$ is moreover big, then $H^i(Y,\inv{\mathcal L})$ vanishes for $i<\op{dim}Y$. 
\end{stelling}

Theorem~\ref{T:MT} also begs the question what (weakly) F-regular type and \qreg{ity} amount to if we drop the $\mathbb Q$-Gorenstein condition. Without the $\mathbb Q$-Gorenstein assumption, one should actually use the notion of strongly F-regular type (which is only conjecturally equivalent with weakly F-regular type), but even then it is no longer clear that this is equivalent with \qreg{ity} (one direction holds by Proposition~\ref{P:qreg} below). In \cite[\S4.6]{HW}, the authors propose \name{Nakayama}'s notion  of \emph{admissible singularities} (\cite{Nak}) as a candidate for an equivalent condition to strongly F-regular type. They point out that  an affine cone of a  smooth Fano variety has in general only admissible singularities (although its anti-canonical cone has log-terminal singularities). The fact that any such cone is \qreg\ (see Remark~\ref{R:all}) corroborates hence their claim.

\subsection{Remark on Kodaira Vanishing}
Note that \name{Hara}'s proof of implication (\ref{i:lt})$\implies$(\ref{i:lt}') in Theorem~\ref{T:MT} relies heavily on Kodaira Vanishing (in fact, on Akizuki-Kodaira-Nakano Vanishing). Therefore, it is of interest to see which of the results in this paper do not make use of Kodaira Vanishing. If we let $S$ be regular in Theorem~\ref{T:ltpure}, then we do not need the implication (\ref{i:lt})$\implies$(\ref{i:lt}') and hence no Vanishing Theorem is used (see Remark~\ref{R:reg} below). Similarly, our proof of Theorem~\ref{T:KVvan} uses only elementary results from cohomology theory and hence does not rely on Kodaira Vanishing. Nonetheless,   in order to prove that smooth Fano varieties are globally \qreg, and hence to obtain Theorem~\ref{T:Fanoquot}, we do need \name{Hara}'s result and hence  Kodaira Vanishing.

\subsection{Remark on the base field}\label{R:bf}
To make the exposition more transparent, I have only dealt in the text with the case that the base field is $\mathbb C$. However, the results extend to arbitrary base fields of \ch\ zero by the following two observations. First, any uncountable \acf\ of \ch\ zero is the ultraproduct of (algebraically closed) fields of positive \ch\ by the Lefschetz Principle (see for instance \cite[Remark 2.5]{SchNSTC}) and this is the only property we used of $\mathbb C$ (cf.\ \eqref{eq:LP} below). Second, since all properties admit faithfully flat descent, we can always make a base change to an uncountable \acf.

\subsection*{Acknowledgement}
I am grateful to Karen Smith for drawing my attention to the question of Vanishing Theorems for GIT quotients of Fano varieties. Without her continuous help and encouragement, this paper would not have been possible.

\section{Transfer and \Sr{s}}

In this section, I will briefly discuss an alternative construction of the usual 'reduction modulo $p$' construction from algebraic geometry. The advantage is that we can work just with schemes of finite type over   fields, that is to say, there is no need to work with relative versions. The one drawback is that we need the base field to be uncountable and algebraically closed (or more generally, an ultraproduct of fields). However, as explained in Remark~\ref{R:bf}, this is not too much of a constraint. Moreover, in order to simplify the exposition, I will in the sequel only discuss the case that the base field is $\mathbb C$.

For generalities on ultraproducts, including \los, see \cite[\S 2]{SchNSTC}. Recall that an ultraproduct  of rings $\seq Cp$ is a certain homomorphic image of the direct product of the $\seq Cp$. This ultraproduct will be denoted by $\up p\seq Cp$, or simply by $\ul C$, and similarly, the image of a sequence $(\seq ap\mid p)$ in $\ul C$ will be denoted by $\up p\seq ap$, or simply by $\ul a$. Any choice of sequence of elements $\seq ap$ whose ultraproduct is equal to $\ul a$ will be called an \emph\sr\ of $\ul a$ (note that we are using the term more loosely than in  \cite{SchNSTC}, where we reserved the notion of \sr\ only for \emph{standard} $\ul a$). The key ingredient for transfer between zero and positive \ch\  is the following fundamental isomorphism
	\begin{equation}\label{eq:LP}
		\mathbb C\iso \up p \ac Fp,
	\end{equation}
where $\ac Fp$ denotes the algebraic closure of the $p$-element field. I will refer to \eqref{eq:LP} as the \emph{Lefschetz Principle} for \acf{s}; see \cite[Theorem 2.4]{SchNSTC} or \cite[Fact 4.2]{SchBS} for proofs.

\subsection{Affine Algebras}\label{s:affalg}
Let me briefly recall from \cite{SchNSTC} the construction of an \sr\ of a finitely generated (for short, an \emph{affine}) $\mathbb C$-algebra $C$. For a fixed tuple of variables  $X$, let $\ul A$ be the ultraproduct of the $\seq Ap:=\pol{\ac Fp}X$. We call $\ul A$ the \emph{non-standard hull} of $A:=\pol{\mathbb C}X$ and $\seq Ap$ an \emph\sr\ of $A$. By \cite{vdD79}, the canonical \homo\ $A\to \ul A$ is faithfully flat (see also \cite[Theorem 1.7]{SvdD} or \cite[A.2]{SchBArt}).  For an arbitrary affine $\mathbb C$-algebra $C$, say of the form  $A/I$, we let 
	\begin{equation}\label{eq:nshull}
	\ul C:=\ul A/I\ul A=C\tensor_A\ul A
	\end{equation}
and call it the \emph{non-standard hull} of $C$. One shows that $\ul C$ is the ultraproduct of affine $\ac Fp$-algebras $\seq Cp$. Any such choice of $\seq Cp$ is called an \emph\sr\ of $C$. There are two ways to construct these: either one observes that $I\ul A$ is the ultraproduct of ideals $\seq Ip$ in $\seq Ap$ (we call $\seq Ip$ an \emph\sr\ of $I$; see \cite[\S3]{SchNSTC}) and takes $\seq Cp$ to be $\seq Ap/\seq Ip$, or alternatively, one takes a model of $C$,  looks at its reductions modulo $p$ and takes a suitably chosen base change over $\ac Fp$ (see Proposition~\ref{P:model} and \S\ref{R:sr} below). By \cite[3.4]{SchNSTC}, the non-standard hull $\ul C$ is independent (up to an isomorphism of $\mathbb C$-algebras) of the presentation of $C$ as a homomorphic image of a polynomial ring. Consequently, if $\seq Cp'$ are obtained from $C$ by the above process starting from a different presentation of $C$, then $\seq Cp\iso \seq Cp'$ (as $\ac Fp$-algebras), for almost all $p$. The non-standard hull does depend though on the choice of ultrafilter and on the choice of the isomorphism~\eqref{eq:LP}. 

We can similarly define the non-standard hull  for a local ring $R$ essentially of finite type over $\mathbb C$ (a \emph{local affine} $\mathbb C$-algebra, for short). Suppose $R$ is of the form $C_\pr$ with $C$ an affine $\mathbb C$-algebra and $\pr$ a prime ideal of $C$. It follows from \cite[Corollary 4.2]{SchNSTC} that $\pr\ul C$ is prime. We define the \emph{non-standard hull} of $R$ to be the localization
	\begin{equation*}
	\ul  R:=(\ul C)_{\pr\ul R}.
	\end{equation*}
This is again independent from the choice of presentation $R=C_\pr$. As explained above, $\pr\ul C$ is an ultraproduct of ideals $\seq\pr p$ in an \sr\ $\seq Cp$ of $C$. By \los, almost all $\seq \pr p$ are prime. We call the localization $\seq Rp:=(\seq Cp)_{\seq\pr p}$ (for those $p$ for which it makes sense), an \emph\sr\ of $R$. It follows that the ultraproduct of the $\seq Rp$ is $\ul R$.

\subsection{Homomorphisms}\label{s:hom}
Let   $\varphi\colon C\to D$ be a (local) \homo\ of finite type between (local) affine $\mathbb C$-algebras. This corresponds to a presentation of $D$ as $\pol CX/I$ (or a localization of the latter), for some   finite tuple of variables $X$. Let $\seq Cp$ and $\seq Dp$ be \sr{s} of $C$ and $D$ respectively, where we use the presentation $D=\pol CX/I$ to construct the $\seq Dp$. This shows that almost all $\seq Dp$ are $\seq Cp$-algebras. The corresponding ring \homo\ $\seq\varphi p\colon \seq Cp\to \seq Dp$ is called  an \emph\sr\  of $\varphi$. The ultraproduct of the $\seq\varphi p$ is a \homo\ $\ul\varphi\colon \ul C\to \ul D$, called the \emph{non-standard hull} of $\varphi$, where $\ul C$ and $\ul D$ are the non-standard hulls of $C$ and $D$ respectively. We have a commutative diagram
	\commdiagram [nsh] C\varphi D {} {} {\ul C} {\ul\varphi} {\ul D.}
Note that if we choose a polynomial ring $A$ of which both $C$ and $D$ are homomorphic images, then $\ul C\iso C\tensor_A\ul A$ and $\ul D\iso D\tensor_A\ul A$ by \eqref{eq:nshull}, and $\ul\varphi$ is just the base change of $\varphi$ over $\ul A$.

\subsection{Affine Schemes}\label{s:affsch}
If $X$ is an affine scheme of finite type over $\mathbb C$, say of the form $\op{Spec} C$ with $C$ an affine $\mathbb C$-algebra, then we call $\seq Xp:=\op{Spec}\seq Cp$ an \emph\sr\ of $X$, for any choice of \sr\ $\seq Cp$ of $C$. One has to be careful however: it is not true that the ultraproduct $\up p\seq Xp$ of the $\seq Xp$ is equal to $\op{Spec}\ul C$. In fact, $\up p\seq Xp$ is the subset of $\op{Spec}\ul C$ consisting of all prime ideals of the form $(\ul I:\ul a)$ for $\ul I$ a finitely generated ideal in $\ul C$ and $\ul a$ an element of $\ul C$ (and hence in general is no longer a scheme). Instead, we call $\ul X:=\op{Spec}(\ul C)$ the \emph{non-standard hull} of $X$. We have a faithfully flat canonical morphism $\ul X\to X$. Since $\ul X$ is no  longer a Noetherian scheme, it is more prudent to reason on its \sr{s} $\seq Xp$ instead, and that is the course we will take in this paper.

\subsection{Affine Morphisms}\label{s:affmor}
Let $f\colon Y\to X$ be a  morphism of finite type between the affine schemes $Y=\op{Spec}D$ and $X=\op{Spec} C$ of finite type over $\mathbb C$. This induces a $\mathbb C$-algebra \homo\ $\varphi\colon C\to D$. Let $\seq \varphi p\colon \seq Cp\to \seq Dp$ be an \sr\ of $\varphi$ (as in \S\ref{s:hom}) and let $\seq fp\colon \seq Yp\to \seq Xp$ be the corresponding morphism between the \sr{s} $\seq Yp:=\op{Spec}\seq Dp$ and $\seq Xp:=\op{Spec}\seq Cp$. We call $\seq fp$ an \emph\sr\ of $f$.  It follows from the corresponding transfer for affine algebras (see \cite[\S4]{SchNSTC}) that if $f$ is an (open, closed, locally closed) immersion (respectively, injective, surjective, an isomorphism, flat, faithfully flat), then so are almost all $\seq fp$. We leave the details to the reader.

\subsection{Modules}\label{s:mod}
Let $\mathcal F$ be a coherent $\loc_X$-module. Any such module is of the form $\widetilde M$ with $M$ a finitely generated $C$-module. Write $M$ as the cokernel of a matrix $\Gamma$ over $C$, that is to say, given by an exact sequence
	\begin{equation*}
	C^a\map\Gamma  C^b \to M\to 0.
	\end{equation*}
Let $\seq \Gamma p$ be an \sr\ of $\Gamma$ (that is to say, the $\seq \Gamma p$ are $(a\times b)$-matrices over $\seq Cp$ with ultraproduct equal to  $\Gamma$) and let $\seq Mp$ be the cokernel of $\seq\Gamma p$. We call $\seq Mp$ an \emph\sr\ of $M$ and we call their ultraproduct $\ul M$ the \emph{non-standard hull} of $M$. Again one shows that $\ul M$ does not depend on the choice of matrix $\Gamma$; in fact, we have an isomorphism
	\begin{equation}\label{eq:uliso}
	\ul M\iso M\tensor_C\ul C.
	\end{equation}
The $\loc_{\seq Xp}$-module $\seq{\mathcal F}p:=\widetilde{\seq Mp}$ associated to $\seq Mp$ is called an \emph\sr\ of $\mathcal F$.

\subsection{Schemes}\label{s:sch}
Let $X$ be a scheme of finite type over $\mathbb C$. Let $U_i$ be a finite covering of $X$ by affine open subsets. For each $i$, let $\seq{U_i}p$ be an \sr\ of $U_i$. I claim  that for almost all $p$, the $\seq{U_i}p$ glue together into a scheme $\seq Xp$ of finite type over $\ac Fp$, and, for any other choice of open affine covering $\{U_i'\}$ of $X$, if the resulting glued schemes are denoted $\seq Xp'$, then $\seq Xp\iso \seq Xp'$, for almost all $p$. This justifies calling the $\seq Xp$ an \emph\sr\ of $X$. The proof of the claim is not hard, but is a little tedious, in that we have to check that the whole construction of glueing schemes is constructive and hence passes by \los\ through ultraproducts. Here is a rough sketch: for each pair $i<j$, we have an isomorphism 
	\begin{equation*}
	\varphi_{ij}\colon \restrict{\loc_{U_i}}{U_i\cap U_j} \iso \restrict{\loc_{U_j}}{U_i\cap U_j}.
	\end{equation*}
Taking \sr{s} $\seq{\varphi_{ij}}p$ of these maps as described in \S\ref{s:affmor}, we get from \los\ that $\seq{\varphi_{ij}}p$ defines an isomorphism
	\begin{equation*}
	\restrict{\loc_{\seq{U_i}p}}{\seq{U_i}p\cap \seq{U_j}p} \iso \restrict{\loc_{\seq{U_j}p}}{\seq{U_i}p\cap \seq{U_j}p}
	\end{equation*}
for almost all $p$. Hence the $\seq{U_i}p$ glue together to get a scheme $\seq Xp$. If we start from a different open affine covering $\{U_i'\}$, then to see that the resulting schemes $\seq Xp'$ agree for almost all $p$, reason on a common refinement of these two coverings.

Similarly, if $C_i$ is the affine coordinate ring of $U_i$, then the $\op{Spec}(\ul{C_i})$ glue together and the resulting scheme $\ul X$ will be called the \emph{non-standard hull} of $X$. In particular, the canonical morphism $\ul X\to X$ is faithfully flat (since it is so locally).

\subsection{Morphisms}\label{s:mor}
Let $f\colon Y\to X$ be a morphism of finite type between schemes of finite type over $\mathbb C$. Let $\seq Xp$ and $\seq Yp$ be \sr{s} of $X$ and $Y$ respectively. Choose finite affine open coverings $\mathfrak U$ and $\mathfrak V$ of respectively $X$ and $Y$, such that $\mathfrak V$ is a refinement of $\inverse f{\mathfrak U}$. In other words, for each $V\in\mathfrak V$, we can find $U\in\mathfrak U$, such that $f(V)\subset U$. Let us write $\restrict fV$ for the restriction $V\to U$ induced by $f$. Choose \sr{s} $\seq{\mathfrak U}p$, $\seq{\mathfrak V}p$ and $\seq{(\restrict fV)}p$ of $\mathfrak U$, $\mathfrak V$ and the affine morphisms $\restrict fV$ respectively (use \S\ref{s:affmor} for the latter). It follows that for any two opens $V,V'\in\mathfrak V$, the morphisms $\seq{(\restrict fV)}p$ and $\seq{(\restrict f{V'})}p$ agree on the intersection $\seq Vp\cap\seq Vp'$, for almost all $p$, and therefore determine a morphism $\seq fp\colon \seq Yp\to \seq Xp$, which we will call an \emph\sr\ of $f$. As for affine morphisms, most algebraic properties descend to the \sr{s} in the sense that $f$ has a certain property (such as being a closed immersion or flat) \iff\ almost all $\seq fp$ have.

\subsection{Coherent Sheaves}\label{s:cohsh}
Let $\mathcal F$ be a coherent $\loc_X$-module. For each $i$, let $\seq{\mathcal G_i}p$ be an \sr\ of the coherent $\loc_{U_i}$-module $\restrict{\mathcal F}{U_i}$ as explained in \S\ref{s:mod} and \S\ref{s:sch} and with the notations therein. Again one easily checks that these $\seq{\mathcal G_i}p$ glue together to give rise to a coherent $\loc_{\seq Xp}$-module $\seq {\mathcal F}p$, which we therefore call an \emph\sr\ of $\mathcal F$, and, moreover, the construction does not depend on the choice of open affine covering.

If $\sheaf$ is a coherent sheaf of ideals on $X$ with \sr\ $\seq{\mathcal F}p$,   almost all $\seq{\mathcal F}p$ are sheaves of ideals, and the closed subscheme they determine on $\seq Xp$ is an \sr\ of the closed subscheme determined by $\mathcal F$. More generally, many local properties (such as being invertible, locally free) hold for the sheaf $\sheaf$ \iff\ they hold for almost all of its \sr{s} $\seq\sheaf p$, since they can be checked locally and hence reduce to a similar transfer property for affine algebras discussed at large in \cite{SchNSTC}.

\subsection{Graded Rings and Modules}\label{s:grad}
Recall   that a ring $S$ is called ($\zet$-)\emph{graded} if it can be written as a direct sum $\oplus_{j\in\zet} \grad Sj$, where each $\grad Sj$ is an additive subgroup of $S$ (called the \emph{$j$-th homogeneous piece} of $S$), with the property  that $\grad Si\cdot \grad Sj\subset \grad S{i+j}$, for all $i,j\in\zet$. In particular, each $\grad Si$ is an $\grad S0$-module. If all $\grad Sj$ are zero for $j<0$, we call $S$ \emph{positively graded}. An $S$-module  $M$ is called \emph{graded} if it admits a decomposition $\oplus_{j\in\zet} \grad Mj$, where each $\grad Mj$ is an additive subgroup of $M$ (called the \emph{$j$-th homogeneous piece} of $M$), with the property  that $\grad Si\cdot \grad Mj\subset \grad M{i+j}$, for all $i,j\in\zet$. We will write $M(m)$ for the \emph{$m$-th twist} of $M$, that is to say, for the graded $S$-module for which $\grad{M(m)}j=\grad M{m+j}$.

Let $S$ be a graded affine $\mathbb C$-algebra. Let $\seq Sp$ be an \sr\ of $S$ and $\ul S$ its non-standard hull. Our goal is to show that almost all $\seq Sp$ are graded. Let $x_i$ be homogeneous algebra generators of $S$ over $\mathbb C$, say of degree $d_i$. Put $A:=\pol{\mathbb C}X$ and let $\varphi\colon A\to S$ be given by $X_i\mapsto x_i$. We make $A$ into a graded ring be giving $X_i$ weight $d_i$, that is to say, $\grad Aj$ is the vector space over $\mathbb C$ generated by all monomials $X_1^{e_1}\cdots X_n^{e_n}$ such that $d_1e_1+\dots+d_ne_n=j$. Hence the kernel $I$ of $A\to S$ is generated by homogeneous polynomials in this new grading. Give each $\seq Ap$ the same grading as $A$ (using the weights $d_i$) and let $\seq Ip$ be an \sr\ of $I$. It follows from \los\ that almost all $\seq Ip$ are generated by homogeneous elements. Since $\seq Sp\iso \seq Ap/\seq Ip$ for almost all $p$, we proved that almost all \sr{s} are graded.   Moreover, if $S$ is positively graded, then so are almost all $\seq Sp$. 

However, the non-standard hull $\ul S$ is no longer a graded ring. Nonetheless, for each non-standard integer $j$ (that is to say, any element $j:=\up p\seq jp$ in the ultrapower $\ul\zet$), we can define the \emph{$j$-th homogeneous piece} $\grad{\ul S}j$ of $\ul S$ as the ultraproduct of the $\grad{\seq Sp}{\seq jp}$. It follows that each $\grad{\ul S}j$ is a direct summand of $\ul S$ (and in fact, $\gr {\ul S}:=\oplus_{j\in\ul\zet}\grad {\ul S}j$ is a direct summand of $\ul S$), and $\grad{\ul S}i\cdot \grad{\ul S}j\subset \grad{\ul S}{i+j}$, for all $i,j\in\ul\zet$ (so that $\gr{\ul S}$ is a $\ul\zet$-graded ring). If $j$ is a standard integer (that is to say, $j\in\zet$,  whence $\seq jp=j$ for almost all $p$), the embedding $S\subset \ul S$ induces an embedding
	\begin{equation}\label{eq:nshgrad}
	\grad Sj\subset \grad {\ul S}j.
	\end{equation}
Note that this is not necessarily an isomorphism. For instance, if $S=\pol{\mathbb C}{X,Y,1/X}$ with $X$ and $Y$ having weight $1$ (and $1/X$ weight $-1$), then $\grad S0\iso\pol{\mathbb C}{Y/X}$ whereas $\grad{\ul S}0$ contains for instance  the ultraproduct of the $Y^p/X^p$.

Let $M$ be a finitely generated graded $S$-module. Let $\seq Mp$ be an \sr\ of $M$ and   $\ul M$ its non-standard hull. By the same argument as above, almost all $\seq Mp$ are graded $\seq Sp$-modules. We define similarly the \emph{$j$-th homogeneous piece} $\grad{\ul M}j$ of $\ul M$ as the ultraproduct of the $\grad{\seq Mp}{\seq jp}$. It follows that $\grad {\ul S}i\cdot \grad{\ul M}j\subset \grad{\ul M}{i+j}$, for each $i,j\in\ul\zet$, and $\grad Mj\subset \grad {\ul M}j$ for each standard $j$. 

If $M\to N$ is a degree preserving morphism of finitely generated graded $S$-modules (so that $\grad Mj$ maps inside $\grad Nj$, for all $j$), then the same is true for almost all \sr{s} $\seq Mp\to \seq Np$. Hence the base change $\ul M\to \ul N$ sends $\grad{\ul M}j$ inside $\grad{\ul N}j$, for each $j\in\ul\zet$.

\subsection{Projective Schemes}\label{s:projsch}
Suppose that $X=\op{Proj}S$ is a projective scheme, with $S$ an affine, positively graded $\mathbb C$-algebra. Let $\seq Sp$ be an \sr\ of $S$. Then $\seq Sp$ is an affine, positively graded $\ac Fp$-algebra by \S\ref{s:grad}, and $\seq Xp\iso\op{Proj}\seq Sp$,  for almost all $p$. Indeed,  this is clear for $S=\pol{\mathbb C}{X_0,\dots,X_n}$ (so that $X=\mathbb P_{\mathbb C}^n$), and the general case follows from this since any projective scheme of finite type over $\mathbb C$ is a closed subscheme of some $\mathbb P_{\mathbb C}^n$.

\subsection{Polarizations}\label{s:polar}
Let $X$ be a projective variety and $\mathcal P$ an  ample line bundle on $X$ (we will call $\mathcal P$ a \emph{polarization} and study this situation in more detail in \S\ref{s:pol}). Let $\seq Xp$ and $\seq{\mathcal P}p$ be \sr{s} of $X$ and $\mathcal P$ respectively. I claim that almost all  $\seq{\mathcal P}p$ are polarizations. That almost all are invertible is clear from the discussion in \S\ref{s:cohsh}. Suppose first that $\mathcal P$ is very ample. Hence there is an embedding $f\colon X\to \mathbb P^N_{\mathbb C}$ such that $\mathcal P\iso f^*\loc(1)$. From the discussion in \S\ref{s:mor} and \S\ref{s:grad}, the \sr\ $\seq fp\colon \seq Xp\to \mathbb P^N_{\ac Fp}$ is an embedding and $\seq{\mathcal P}p\iso \seq fp^*\loc(1)$ for almost all $p$, showing that almost all $\seq{\mathcal P}p$ are very ample. If $\mathcal P$ is just ample, then $\mathcal P^m$ is very ample for some $m>0$ by \cite[II. Theorem 7.6]{Hart}. Hence by our previous argument, almost all $\seq{\mathcal P}p^m$ are very ample. By another application of \cite[II. Theorem 7.6]{Hart}, almost all $\seq{\mathcal P}p$ are ample. 

Presumably the converse also holds, but this requires a finer study of the dependence of the  exponent $m$ on the ample sheaf: it should only depend on the degree complexity of the sheaf (that is to say, on the maximum of the degrees of the polynomials needed in describing the sheaf).

\subsection{Complexes}\label{s:comp}
Let $\cc$ be an arbitrary bounded complex in which each term $\mathcal C^m$ is a finitely generated module over an affine $\mathbb C$-algebra. Using \S\ref{s:hom} and \S\ref{s:affsch}, we can choose  an \sr\ for each term and each \homo\ in this complex. Let $\seq{(\cc)} p$ denote the corresponding object. By \los, almost all $\seq{(\cc)} p$ are complexes. This justifies calling $\seq{(\cc)} p$  an \emph\sr\ of $\cc$. Let $A$ be a polynomial ring over $\mathbb C$ such that each $\mathcal C^m$ is an $A$-module. It follows from \eqref{eq:uliso} that we have an isomorphism of complexes
	\begin{equation}\label{eq:ulcom}
	\cc\tensor_A\ul A\iso \up p \seq{(\cc)}p.
	\end{equation}
Since taking cohomology consists of taking kernels, images  and quotients, each of which commutes with ultraproducts, taking cohomology also commutes with ultraproducts. Applying this to \eqref{eq:ulcom}, we get for each $i$, an  isomorphism 
	\begin{equation}\label{eq:ulcoh}
	H^i(\cc) \tensor_A\ul A\iso H^i(\cc\tensor_A\ul A)\iso \up p H^i(\seq{(\cc)}p)
	\end{equation}
where we used that $A\to \ul A$ is faithfully flat in the first isomorphism.

Our next goal is to show that an \sr\ of the cohomology of a coherent $\loc_X$-module $\mathcal F$ is obtained by taking the cohomology of its \sr{s}. In order to prove this, we use \Cech\ cohomology to calculate sheaf cohomology (this will be studied further in  \S\ref{ss:ccs} below).

\subsection{\Cech\ Cohomology}\label{s:cecoh}
Recall that the \emph{\Cech\ complex} $\cc(\mathfrak U;\mathcal F)$ of $\mathcal F$ associated to an open affine covering $\mathfrak U:=\{U_1,\dots,U_s\}$ of $X$ is by definition the complex in which the $m$-th term for $m\geq 1$ is
	\begin{equation*}
	\mathcal C^m(\mathfrak U;\mathcal F):= \bigoplus_{\mathbf i} H^0(U_{\mathbf i}, \mathcal F)
	\end{equation*}
where $\mathbf i$ runs over all $m$-tuples of indices $1\leq i_1<i_2<\dots<i_m\leq s$ and where $U_{\mathbf i}:=U_{i_1}\cap U_{i_2}\cap \dots\cap U_{i_s}$ (see \cite[III.4]{Hart} for more details). Note that $\cc(\mathfrak U;\mathcal F)$ is a  bounded complex of affine $\mathbb C$-algebras. 

\begin{lemma}\label{L:srcomp}
Let $X$ be a scheme of finite type over $\mathbb C$, let $\mathfrak U$ be a finite affine open covering of $X$ and let $\mathcal F$ be a coherent $\loc_X$-module. If $\seq Xp$, $\seq{\mathfrak U}p$ and $\seq {\mathcal F}p$ are \sr{s} of $X$, $\mathfrak U$ and $\mathcal F$ respectively, then the complexes $\cc(\seq{\mathfrak U}p; \seq{\mathcal F}p)$ are  an \sr\ of the complex $\cc(\mathfrak U;\mathcal F)$.
\end{lemma}
\begin{proof}
Since an \sr\ of $\mathfrak U$ is obtained by choosing an \sr\ for each affine open in it, we get from \los\ that $\seq{\mathfrak U}p$ is an open covering of $\seq Xp$ for almost all $p$. Moreover, if $U$ is an affine open with \sr\ $\seq Up$, then $H^0(\seq Up,\seq{\mathcal F}p)$ is an \sr\ of $H^0(U,\mathcal F)$. The assertion readily follows from these observations.
\end{proof}

If $X$ is separated and of finite type over $\mathbb C$ and if  $\mathcal F$ is a coherent $\loc_X$-module,  then the cohomology modules $H^i(X,\mathcal F)$ can be calculated as the cohomology of the \Cech\ complex $\cc(\mathfrak U;\mathcal F)$, for any choice of finite  open affine covering $\mathfrak U$ (\cite[Theorem 4.5]{Hart}). More precisely, 
	\begin{equation}\label{eq:ccoh}
	H^i(X,\mathcal F)\iso H^{i+1}(\cc(\mathfrak U;\mathcal F))
	\end{equation}
(some authors start numbering the \Cech\ complex from zero, so that there is no shift in the superscripts needed). If $\mathfrak U$ consists of affine opens $\op{Spec}C_i$, we can choose   a polynomial ring $A$ over $\mathbb C$ so that every $C_i$ is a homomorphic of $A$. It follows that each module in $\cc(\mathfrak U;\mathcal F)$ is a finitely generated $A$-module, and hence so is each $H^i(X,\mathcal F)$. 

If $X$ is moreover projective, then each $H^i(X,\mathcal F)$ is a finite dimensional vector space over $\mathbb C$ and its dimension will be denoted by $h^i(X,\mathcal F)$.

\begin{theorem}\label{T:nscoh}
Let $X$ be a separated scheme of finite type over $\mathbb C$  and let $\mathcal F$ be a coherent $\loc_X$-module. Let $\seq Xp$ and $\seq {\mathcal F}p$ be \sr{s} of $X$ and $\mathcal F$ respectively. For an appropriate choice of a  polynomial ring $A$ over $\mathbb C$ and for each $i$, we have an isomorphism
	\begin{equation}\label{}
	H^i(X,\mathcal F)\tensor_A\ul A\iso \up p H^i(\seq Xp,\seq{\mathcal F}p).
	\end{equation} 

In particular, if $X$ is moreover projective, then $h^i(X,\mathcal F)$ is equal to $h^i(\seq Xp,\seq{\mathcal F}p)$ for almost all $p$.
\end{theorem}
\begin{proof}
By Lemma~\ref{L:srcomp}, we have for almost all $p$ an isomorphism of complexes
	\begin{equation*}
	\seq{(\cc(\mathfrak U;\mathcal F))}p \iso \cc(\seq{\mathfrak U}p; \seq{\mathcal F}p).
	\end{equation*}
where the left hand side is some \sr\ of $\cc(\mathfrak U;\mathcal F)$. The first assertion now follows from \eqref{eq:ulcoh} and \eqref{eq:ccoh}. The last assertion follows from the first, by taking lengths of both sides and using \cite[Proposition 1.5]{SchEC}.
\end{proof}

In particular, if $H^i(X,\mathcal F)$ vanishes for some $i$, then so will almost all $H^i(\seq Xp,\seq{\mathcal F}p)$. More precisely, for a fixed choice of \sr, let $\varpi_i(\mathcal F)$ be the collection of prime numbers $p$ for which $H^i(\seq Xp,\seq{\mathcal F}p)$ vanishes. By the above result, $\varpi_i(\mathcal F)$ belongs to the ultrafilter \iff\ $H^i(X,\mathcal F)=0$. However, if we have an infinite collection of coherent sheaves $\mathcal F_n$ with zero $i$-th cohomology, the intersection of all $\varpi_i(\mathcal F_n)$ will in general no longer belong to the ultrafilter, and therefore can very well be empty. The next result shows that by imposing some further algebraic relations among the $\mathcal F_n$, the intersection remains in the ultrafilter.

\begin{corollary}\label{C:van}
Let $X$ be a projective scheme of finite type over $\mathbb C$. Let $\mathcal L$ be an invertible $\loc_X$-module and let $\mathcal E$ be a locally free $\loc_X$-module. Let $\seq Xp$, $\seq{\mathcal L}p$  and $\seq{\mathcal E}p$ be \sr{s} of $X$, $\mathcal L$  and $\mathcal E$ respectively. If for some $i$ and some $n_0$, we have that $H^i(X,\mathcal E\tensor\mathcal L^n)=0$ for all $n\geq n_0$, then for almost all $p$, we have, for all $n\geq n_0$, that $H^i(\seq Xp,\seq{\mathcal E}p\tensor\seq{\mathcal L}p^n)=0$. 
\end{corollary}
\begin{proof}
Let $\mathcal A$ denote the symmetric algebra $\oplus_{n\geq 0} \mathcal L^n$ of $\mathcal L$ and let $\mathcal F:=\mathcal A\tensor\mathcal E\tensor\mathcal L^{n_0}$.  Note that $\mathcal F=\oplus_{n\geq n_0} \mathcal E\tensor\mathcal L^n$, so that our assumption becomes $H^i(X,\mathcal F)=0$. We cannot apply Theorem~\ref{T:nscoh} directly, as $\mathcal F$ is not a coherent $\loc_X$-module.  Let $Y:=\op{\bf Spec}\mathcal A$ be the scheme over $X$  associated to $\mathcal A$ (see \cite[II. Ex. 5.17]{Hart}). Since $\mathcal A$ is a finitely generated sheaf of $\loc_X$-algebras, the morphism $f\colon Y\to X$ is of finite type. Moreover, $f$ is affine (that is to say, $\inverse fU\iso \op{Spec} \mathcal A(U)$ for every affine open $U$ of $X$) and $\mathcal A\iso f_*\loc_Y$. Let $\mathcal G:=f^*(\mathcal E\tensor\mathcal L^{n_0})$, so that $\mathcal G$ is a coherent $\loc_Y$-module.  We have isomorphisms
	\begin{equation*}
	f_*\mathcal G\iso f_*\loc_Y\tensor \mathcal E\tensor\mathcal L^{n_0} \iso \mathcal A\tensor \mathcal E\tensor\mathcal L^{n_0}=\mathcal F,
	\end{equation*}
where the first isomorphism follows from the projection formula (see \cite[II. Ex. 5.1]{Hart}). Therefore, 
	\begin{equation}\label{eq:ycoh}
	H^i(Y,\mathcal G)\iso H^i(X,f_*\mathcal G) = H^i(X,\mathcal F)=0
	\end{equation}
where the first isomorphism holds by \cite[III. Ex. 4.1]{Hart}. 

Let $\seq fp\colon \seq Yp\to\seq Xp$ be an \sr\ of $f$ (as described in \S\ref{s:mor}) and let $\seq{\mathcal F}p$ and $\seq{\mathcal G}p$ be \sr{s} of $\mathcal F$ and $\mathcal G$ respectively. By \los, we have isomorphisms
	\begin{align}\label{}
	\seq{\mathcal G}p&\iso \seq fp^*(\seq{\mathcal E}p\tensor\seq{\mathcal L}p^{n_0})\\
	(\seq fp)_*\seq{\mathcal G}p&\iso \seq{\mathcal F}p\\
	(\seq fp)_*\loc_{Y_p}&\iso \bigoplus_{n\geq 0} \seq{\mathcal L}p^n \\
	\seq{\mathcal F}p&\iso \bigoplus_{n\geq n_0}\seq{\mathcal E}p\tensor \seq{\mathcal L}p^n\label{eq:fp}
	\end{align}
for almost all $p$. Applying Theorem~\ref{T:nscoh} to \eqref{eq:ycoh}, we get that almost all $H^i(\seq Yp,\seq{\mathcal G}p)$ vanish. Hence, by the analogue of \eqref{eq:ycoh} in \ch\ $p$, almost all  $H^i(\seq Xp,\seq{\mathcal F}p)$ vanish.   In view of \eqref{eq:fp}, this proves the assertion.
\end{proof}

Let us conclude this section with showing that this non-standard formalism just introduced is closely related to the usual reduction modulo $p$ (as for instance used in the definition of tight closure in \ch\ zero in \cite{HHZero}).

\subsection{Models}
Let $K$ be a field and $R$ a $K$-affine algebra. With a \emph{model of $R$} (called \emph{descent data} in \cite{HHZero}) we mean a pair $(Z,R_Z)$ consisting of a subring $Z$ of $K$ which is finitely generated over $\zet$ and a $Z$-algebra $R_Z$ essentially of finite type, such that $R\iso R_Z\tensor_Z K$. Oftentimes,  we will think of $R_Z$  as being the model. Clearly, the collection of models $R_Z$ of $R$ forms a direct system whose union is $R$. We say that $R$ is \emph{F-rational type} (respectively, is \emph{weakly F-regular type}, or \emph{strongly F-regular type}), if there  exists a model $(Z,R_Z)$, such that $R_Z/\pr R_Z$ is F-rational (respectively, weakly F-regular or strongly F-regular) for an open set $U$ of maximal ideals $\pr$ of $Z$ (note that $R_Z/\pr R_Z$ has positive \ch). See  \cite{HHZero} or \cite[App. 1]{HuTC} for more details.

 The following was proved in \cite[Proposition 4.10]{SchBCM} for local rings; the general case is proven by the same argument.

\begin{proposition}\label{P:model}
Let $R$ be a $\mathbb C$-affine domain. For each finite subset of $R$, we can find a   model $(Z,R_Z)$ of $R$ containing this subset, and, for almost all $p$, a maximal ideal $\seq \pr p$ of $Z$  and a separable extension  $Z/\seq \pr p\subset\ac Fp$, such that the collection of base changes $R_Z\tensor_Z\ac Fp$ gives  an \sr\ of $R$. Moreover, for any $r\in R_Z$, the collection of images of $r$ under the various \homo{s} $R_Z\to R_Z\tensor_Z\ac Fp$ gives an \sr\ of $r$.
\end{proposition}

\subsection{\Sr{s} as Universal Reductions}\label{R:sr}
Suppose $(R_{Z'},Z')$ is another model of $R$ satisfying the assertion of the previous proposition (so that we have \homo{s} $Z'\to \ac Fp$). Since any two \sr{s} agree almost everywhere as mentioned in \S\ref{s:affalg}, we get that $R_Z\tensor_Z\ac Fp\iso R_{Z'}\tensor_{Z'}\ac Fp$ for almost all $p$.

\section{Log-terminal Singularities}\label{s:lt}

 In \cite{HW}, the authors  show that a $\mathbb Q$-Gorenstein $\mathbb C$-affine algebra has log-terminal singularities if it is  weakly F-regular type  (note that weakly F-regular type and strongly F-regular type are equivalent under the $\mathbb Q$-Gorenstein condition by \cite{MacC}), whereas the converse is proved in \cite{HaRat}. In this section we will define a third condition in terms of ultraproducts of Frobenii and prove its equivalence with the other ones. The advantage of the latter property is that it is easily seen to descend under pure \homo{s} (see Proposition~\ref{P:pure}). We recall some terminology first.

\subsection{$\mathbb Q$-Gorenstein Singularities}\label{D:Qgor}
Recall that a normal scheme $X$ is called \emph{$\mathbb Q$-Gorenstein} if some positive multiple of its canonical divisor $K_X$ is Cartier; the least such positive multiple is called the \emph{index} of $X$. If $f\colon \overline X\to X$ is a resolution of singularities of $X$ and $E_i$ are the irreducible components of the exceptional locus, then the canonical divisor $K_{\overline X}$ is numerically equivalent  to $f^*(K_X)+\sum a_iE_i$ (as $\mathbb Q$-divisors), for some $a_i\in\mathbb Q$ ($a_i$ is called the \emph{discrepancy} of $X$ along $E_i$; see \cite[Definition 2.22]{KM}). If all $a_i>-1$,  we call $X$ \emph{log-terminal} (in case we only have a weak inequality, we call $X$ \emph{log-canonical}).

\subsection{\Qfrobs}\label{s:uf}
Any ring $R$ of \ch\ $p$ is endowed with the \emph{Frobenius} endomorphism $\frob p\colon x\mapsto x^p$, and its powers $\frob q:=\frob p^e$, where $q=p^e$. We can therefore view $R$ as a module over itself via the \homo\ $\frob q$, and to emphasize this, we will use the notation $\bc{\frob q}R$ (a notation borrowed from algebraic geometry; other authors use notations such as $R^{F}$, $R^{\frob q}$ or $\mathstrut^eR$). Similarly, for an arbitrary $R$-module $M$, we will write $\bc{\frob q}M$ for the $R$-module structure on $M$ via $\frob q$ (that is to say, $x\cdot m=x^qm$). It follows that $\bc{\frob q}M\iso M\tensor_R\bc{\frob q}R$.

For each prime number $p$, choose a positive integer $\seq ep$ and let $\pi$ be the non-standard integer given as the ultraproduct of the powers $p^{\seq ep}$. To each such $\pi$, we associate an \emph\qfrob\ in the following way. For each  $\mathbb C$-affine domain  $R$   with non-standard hull $\ul R$, consider the \homo\
	\begin{equation*}
	R\to \ul R\colon x\mapsto x^\pi:=\up p (\seq xp)^{p^{\seq ep}}
	\end{equation*}
where $\seq xp$ is an \sr\  of $x$ (one easily checks that this does not depend on the choice of \sr). We will denote this \qfrob\ by $\uf_\pi$, or simply $\uf$; whenever we want to emphasize the ring $R$ on which it operates, we write $\uf_{\pi;R}$ or simply $\uf_R$. This assignment is functorial, in the sense that for any \homo\ $f\colon R\to S$ of finite type, we have a commutative diagram
	\commdiagram [uf] R fS {\uf_R}{\uf_S} {\ul R}{\ul f}{\ul S.}
 Note that $\uf$ is the restriction to $R$ of the ultraproduct of the $\frob {p^{\seq ep}}$. In particular, if all $\seq ep=1$, then the corresponding \qfrob\ $\ulfrob$ was called the \emph{non-standard Frobenius}   in \cite{SchNSTC}. 

Note that each \qfrob\ induces an $R$-module structure on $\ul R$, which we will denote by $\bc\uf\ul R$ (so that $x\cdot \ul r= \uf(x)\ul r$). It follows that $\bc\uf\ul R$ is the ultraproduct of the $(\frob {p^{\seq ep}})\bc\mathstrut\seq Rp$, with $\seq Rp$ an \sr\ of $R$. If $M$ is a finitely generated $R$-module with non-standard hull $\ul M$, then the $R$-module structure on $\ul M\iso M\tensor_R\ul R$   via the action of $\uf$ on the second factor, will be denoted $\bc\uf M$. It follows that $\bc\uf M$ is isomorphic to $M\tensor_R\bc\uf R$ and hence isomorphic to the ultraproduct of the $(\frob {p^{\seq ep}})\bc\mathstrut\seq Mp$.

\begin{definition}
We say that a $\mathbb C$-affine domain $R$ is \emph\qreg, if for each non-zero $c\in R$, we can find an \qfrob\ $\uf$ such that the $R$-module morphism 
	\begin{equation*}
	 \multf c\uf R\colon R\to \bc\uf\ul R\colon x\mapsto c\uf(x)
	\end{equation*}
is pure. 
\end{definition}

For $M$ an $R$-module, we will write $\multf c\uf M\colon M\to M\tensor_R\bc\uf\ul R$ for the base change of $\multf c\uf R$. Since purity is preserved under localization, one easily verifies that the localization of an \qreg\ ring is again \qreg.

\begin{remark}\label{R:regqreg}
If $R$ is normal, so that purity and cyclical purity are the same by \cite[Theorem 2.6]{HoPure}, then purity of $\multf c\uf R$ is equivalent to the condition that  for every $y\in R$ and every ideal $I$   in $R$, if $c\uf(y)\in\uf(I)\ul R$, then $y\in I$. From this and the fact that any \qfrob\ on a regular local ring is flat (same proof as for \cite[Proposition 6.1]{SchNSTC}), one easily checks that a regular (local) $\mathbb C$-affine domain is \qreg. 
\end{remark}

In \cite{SchBCM}, we called a local $\mathbb C$-affine domain with \sr\ $\seq Rp$ \emph\genreg\ (respectively, \emph\genrat), if each ideal $I$ in $R$ (respectively, some ideal $I$ generated by a system of parameters), is equal to its generic tight closure. Recall from \cite{SchNSTC} that an element $x\in R$ lies in the \emph{generic tight closure} of an ideal $I$, if $\seq xp$ lies in the tight closure of $\seq Ip$, for almost $p$, where $\seq Ip$ and $\seq xp$ are  \sr{s} of $I$ and $x$ respectively. We proved in \cite[Theorem C]{SchBCM} that being \genrat\ is equivalent with having rational singularities. 

\begin{proposition}\label{P:qreg}
Let $R$ be a $\mathbb C$-affine domain. 
\begin{enumerate}
\item\label{i:Fregt} If $R$ is strongly F-regular type, then it is \qreg.
\item\label{i:qreg} If $R$ is \qreg, then every localization is  \genreg.
\end{enumerate}
\end{proposition}
\begin{proof}
To prove \eqref{i:Fregt}, let $c$ be a non-zero element of $R$. Let $(Z,R_Z)$ be a model of $R$ containing $c$. By Proposition~\ref{P:model}, there exists for almost all $p$, a maximal ideal $\seq \pr p$ of $R$ and a separable extension $Z/\seq \pr p\subset \ac Fp$, such that $\seq Rp:=R_Z\tensor_Z\ac Fp$ is an \sr\ of $R$ and such that  for any element $r\in R_Z$, its image in $\seq Rp$ under the base change $\seq\gamma p\colon R_Z\to \seq Rp$ is an \sr\ of $r$. In particular, $c$ is the ultraproduct of the $\seq\gamma p(c)$. By definition, we may choose the model in such way that almost all $\seq Sp:=R_Z/\seq \pr p R_Z$ are strongly F-regular. In particular, we can find powers $q:=p^{\seq ep}$, such that the morphism 
	\begin{equation*}
	 \seq Sp \to \bc{\frob q}(\seq Sp)\colon x\mapsto \seq\gamma p(c)x^q
	\end{equation*}
 is pure. By base change,   the $\seq Rp$-module morphism 
	\begin{equation*}
	 \seq Rp \to \bc{\frob q}(\seq Sp)\tensor\seq Rp\colon x\mapsto \seq\gamma p(c)x^q\tensor 1
	\end{equation*}
 is also pure.  Since $Z/\seq \pr p\subset \ac Fp$ is separable, we get that  $\bc{\frob q}(\seq Sp)\tensor\seq Rp\iso \bc{\frob q}\seq Rp$, showing that the $\seq Rp$-module morphism 
	\begin{equation}\label{eq:cp}
	\seq Rp \to \bc{\frob q}\seq Rp\colon x\mapsto \seq\gamma p(c)x^q
	\end{equation}
 is pure. Let $\uf$ be the \qfrob\ given as the (restriction to $R$ of the) ultraproduct of the $\frob q=\frob p^{\seq ep}$ and let $\ul g$ be the ultraproduct of the morphisms given in \eqref{eq:cp}. It follows that $\ul g(\ul x)=c\uf(\ul x)$, so that the restriction of $\ul g$ to $R$ is precisely $\multf c\uf R$. Moreover,  from the purity of \eqref{eq:cp}, it follows,  using \los, that every finitely generated ideal $J$ of $\ul R$ is  equal to the contraction of its extension under $\ul g$. Since $R\to \ul R$ is faithfully flat, whence cyclically pure, it follows that the restriction of $\ul g$ to $R$, that is to say, $\multf c\uf R$,  is cyclically pure. Since $R$ is in particular normal, $\multf c\uf R$ is pure by \cite[Theorem 2.6]{HoPure}, showing that $R$ is \qreg.

Assume next that $R$ is \qreg. Without loss of generality, we may assume that $R$ is moreover local. Let $I$ be an ideal in $R$ and let $x$ be an element in the generic tight closure of $I$. We need to show that $x\in I$. By \cite[Proposition 8.3]{SchNSTC}, we can choose $c\in R$ such that almost all $\seq cp$ are test elements for $\seq Rp$, where $\seq Rp$ and $\seq cp$ are \sr{s} of $R$ and $c$ respectively. Let $\uf$ be an \qfrob\ such that the $R$-module morphism $\multf c\uf R$ is pure. In particular this implies for every $y\in R$ that
	\begin{equation}\label{eq:pur0}
	\text{if $c\,\uf(y)\in\uf(I)\ul R$, then $y\in I$}.
	\end{equation}
Suppose  $\uf$ is the ultraproduct of the $\frob p^{\seq ep}$. Therefore, \eqref{eq:pur0} translated in terms of an \sr\ $\seq yp$ of an element $y\in R$, becomes the statement 
	\begin{equation}\label{eq:purp}
	\text{if $\seq cp\,\frob p^{\seq ep}(\seq yp)\in\frob p^{\seq ep}(\seq Ip)\seq Rp$, then $\seq yp\in\seq Ip$},
	\end{equation}
for almost all $p$. 

Let $\seq xp$ and $\seq Ip$ be \sr{s} of $x$ and $I$ respectively. By assumption, almost all $\seq xp$ lie in the  tight closure of $\seq Ip$. Since $\seq cp$ is a test element, this means that 
	\begin{equation*}
	\seq cp\,\frob p^N(\seq xp) \in\frob p^N(\seq Ip)\seq Rp,
	\end{equation*}
for all $N$.  With $N=\seq ep$, we get from \eqref{eq:purp} that $\seq xp\in\seq Ip$. Taking ultraproducts,  $x$ lies in $I\ul R$, whence in $I$, by the faithful flatness of $R\to \ul R$.
\end{proof}

\begin{remark}
Note that if $R$ is \qreg\ with \sr\ $\seq Rp$, then it is not necessary the case that almost all $\seq Rp$ are strongly F-regular. Namely,  the set $\Sigma_{y,I}$ of prime numbers $p$ for which \eqref{eq:purp} holds, depends a priori on $y$ and $I$, and therefore, their intersection over all possible $y$ and $I$ might very well be empty.
\end{remark}

The prime \ch\ analogue of the next result was first observed in   \cite{WatF}; we follow the argument given in \cite[Theorem 4.15]{SmVan}.

\begin{proposition}\label{P:Fpure}
Let $R\subset S$ be a finite extension of   local $\mathbb C$-affine domains, \'etale in codimension one. Let $c$ be a non-zero element of $R$ and $\uf$ an \qfrob. If $\multf c\uf R\colon R\to \bc\uf\ul R$ is pure, then so is  $\multf c\uf S\colon S\to \bc\uf\ul S$. 

In particular, if $R$ is \qreg, then so is $S$.
\end{proposition}
\begin{proof}
 Let $R\subset S$ be an arbitrary finite extension  of  $d$-dimensional local $\mathbb C$-affine domains and fix a non-zero element $c$ and an \qfrob\ $\uf$. Let $\mathfrak n$ be the maximal ideal of $S$ and $\omega_S$ its canonical module. I claim that if $R\subset S$ is \'etale, then $S\tensor_R\bc\uf \ul R\iso \bc\uf \ul S$. Assuming the claim, let $R\subset S$ now only be \'etale in codimension one.  It follows from the claim that the supports of the kernel and the cokernel of the base change $S\tensor_R\bc\uf\ul R\to \bc\uf\ul S$ have  codimension at least two. Hence the same is true for the base change
	\begin{equation*}
	\omega_S\tensor_S S\tensor_R\bc\uf{\ul R}\to \omega_S\tensor_S\bc\uf{\ul S}.
	\end{equation*}
Applying the top local cohomology functor $H_{\mathfrak n}^d$, we get, in view of Grothen\-dieck Vanishing and the long exact sequence of local cohomology, an isomorphism
	\begin{equation}\label{eq:top}
	H_{\mathfrak n}^d(\omega_S\tensor_R\bc\uf{\ul R}) \iso H_{\mathfrak n}^d(\omega_S\tensor_S\bc\uf{\ul S}).
	\end{equation}
By Grothendieck duality, $H_{\mathfrak n}^d(\omega_S)$ is the injective hull $E$ of the  residue field of $S$. Taking the base change of $\multf c\uf R$ and $\multf c\uf S$ over $\omega_S$,  and then taking the top local cohomology,  yields the following commutative diagram
	\begin{equation*}
	\CD
	E=H_{\mathfrak n}^d(\omega_S) @>>> E\tensor_R\bc\uf{\ul R} @>>> H_{\mathfrak n}^d(\omega_S\tensor_R\bc\uf{\ul R}) \\
	@|	@VVV	@VV\iso V \\
	E=H_{\mathfrak n}^d(\omega_S) @>>> E\tensor_S\bc\uf{\ul S} @>>> H_{\mathfrak n}^d(\omega_S\tensor_S\bc\uf{\ul S})
	\endCD
	\end{equation*}
where the last  vertical arrow in this diagram is the isomorphism~\eqref{eq:top}. Since by assumption, $\multf c\uf R\colon R\to \bc\uf\ul R$ is pure, so is the base change $\omega_S\to \omega_S\tensor_R\bc \uf{\ul R}$. Since purity is preserved after taking cohomology, the  top composite arrow is injective, and hence so is the lower composite arrow. In particular, its first factor  $E\to E\tensor_S\bc\uf{\ul S}$ is injective. Note that this morphism is still given as  multiplication by $c$, and hence is equal to the base change $\multf c\uf E$ of $\multf c\uf S$. By   \cite[Lemma 2.1(e)]{HHbigCM2}, the injectivity of $\multf c\uf E=\multf c\uf S\tensor E$  then implies  that $\multf c\uf S$ is pure, as we set out to prove.

To prove the claim, observe that if $R\to S$ is \'etale with \sr\ $\seq Rp\to\seq Sp$, then almost all  of these are \'etale. Indeed, by   \cite[Corollary 3.16]{Milne}, we can write $S$ as $\pol RX/I$, with $X=\rij Xn$ and $I=\rij fn\pol RX$, such that the Jacobian $\op J(f_1,\dots,f_n)$   is a unit in $R$, and by \los, this property is preserved for almost all \sr{s}. In general, if $C \to D$ is an \'etale extension of \ch\ $p$ domains, then we have an isomorphism $\bc{\frob q}C\tensor_CD\iso \bc{\frob q}D$ (see for instance \cite[p. 50]{HHTC} or the proof of \cite[Theorem 4.15]{SmVan}). Applied to the current situation, we get that $\seq Sp\tensor_{\seq Rp}\bc{\frob q}\seq Rp\iso\bc{\frob q}\seq Sp$, for $q$ any power of $p$ (\cite[p. 50]{HHTC}). Therefore,  after taking ultraproducts, we obtain the required isomorphism $S\tensor_R\bc\uf\ul R\iso \bc\uf\ul S$ (note that $\ul S\iso \ul R\tensor_RS$ since $R\to S$ is finite).

To prove the last assertion, we have to show that we can find for each non-zero $c\in S$ an \qfrob\ $\uf$ such that $\multf c\uf S$  is pure. However, if we can do this for some non-zero multiple of $c$, then we can also do this for $c$, and hence, since $S$ is finite over $R$,  we may assume without loss of generality that $c\in R$. Since $R$ is \qreg, we can find therefore an \qfrob\ $\uf$ such that $\multf c\uf R$ is pure, and hence by the first assertion, so is then $\multf c\uf S$, proving that $S$ is \qreg. 
\end{proof}

\subsection{Proof of Theorem~\ref{T:MT}}\label{s:pfMT}
The equivalence of (\ref{i:lt}) and (\ref{i:lt}') is proven by  Hara  in \cite[Theorem 5.2]{HaRat}.  Proposition~\ref{P:qreg} proves (\ref{i:lt}')$\implies$(\ref{i:uFreg}). Hence remains to prove (\ref{i:uFreg})$\implies$(\ref{i:lt}). 

To this end,  assume $R$ is \qreg. Recall the construction of the canonical cover of $R$ due to \name{Kawamata}. Let $r$ be the index of $R$, that is to say, the least $r$ such  that $\loc_X(rK_X)\iso \loc_X$, where $X=\op{Spec}R$ and $K_X$ the canonical divisor of $X$. This isomorphism induces an $R$-algebra structure on
	\begin{equation*}
	\tilde R:= H^0(X,\loc_X \oplus \loc_X(K_X)\oplus\dots\oplus \loc_X((r-1)K_X)),
	\end{equation*}
which is called the \emph{canonical cover} of $R$; see \cite{Kaw}.  Since $R\to \tilde R$ is \'etale in codimension one (see for instance \cite[4.12]{SmVan}),  we get from Proposition~\ref{P:Fpure} that $\tilde R$ is \qreg. Hence $\tilde R$ is \genreg, by    Proposition~\ref{P:qreg}. In particular, $\tilde R$ is \genrat, whence has rational singularities, by \cite[Theorem 6.2]{SchRatSing}. By \cite[Theorem 1.7]{Kaw}, this in turn implies that $R$ has log-terminal singularities. \qed

\begin{remark}
Note that without relying on Hara's result (which uses Kodaira Vanishing), we proved the implications (\ref{i:lt}')$\implies$\eqref{i:uFreg}$\implies$\eqref{i:lt}, recovering the result of Smith in \cite{SmFrat,SmVan}. 
\end{remark}

\begin{remark}\label{R:genreg}
There are at least eight more conditions which are expected to be equivalent with the ones in Theorem~\ref{T:MT} for a local $\mathbb Q$-Gorenstein  $\mathbb C$-affine domain $R$, namely 
\begin{enumerate}
\addtocounter{enumi}{2}
\item\label{i:genreg} $R$ is \genreg;
\item[(\ref{i:genreg}')]   $R$ is \genregglob\ (that is to say, every localization of $R$ is \genreg); 
\item\label{i:Freg}  $R$ is weakly F-regular (that is to say, every ideal is equal to its tight closure in the sense of \cite{HHZero});
\item[(\ref{i:Freg}')]   $R$ is F-regular (that is to say, every localization of $R$ is weakly F-regular); 
\item\label{i:Breg}  $R$ is weakly $\mathcal B$-regular (that is to say, $R\to \BCM R$ is cyclically pure; see \S\ref{s:vantor} below); 
\item[(\ref{i:Breg}')]   $R$ is $\mathcal B$-regular (that is to say, every localization of $R$ is weakly $\mathcal B$-regular); 
\item\label{i:nsreg} $R$ is \nsreg\ (that is to say, every ideal is equal to its non-standard tight closure in the sense of \cite{SchNSTC}; see \cite[\S3.10]{SchHR});
\item[(\ref{i:nsreg}')]   $R$ is \nsregglob\ (that is to say, every localization of $R$ is \nsreg); 
\end{enumerate}
The implications \eqref{i:nsreg}$\implies$\eqref{i:genreg}$\implies$\eqref{i:Freg} and \eqref{i:genreg}$\implies$\eqref{i:Breg} follow respectively from the facts that  generic tight closure is contained in non-standard tight closure \cite[Theorem 10.4]{SchNSTC}, that classical tight closure is contained in generic tight closure  \cite[Theorem 8.4]{SchNSTC}, and that $\mathcal B$-closure is contained in generic tight closure \cite[Corollary 4.1]{SchBCM}. We have similar implications among the accented conditions, and, of course,  the accented conditions trivially imply their weak counterparts. In fact, we have even an implication \eqref{i:genreg}$\implies$(\ref{i:Breg}') by \cite[Corollary 4.4]{SchBCM}. Finally, Proposition~\ref{P:qreg} proves that \eqref{i:uFreg}$\implies$(\ref{i:genreg}').

Conjecturally, weakly F-regular is the same as weakly F-regular type, so that therefore all (weak) conditions \eqref{i:lt}--\eqref{i:Freg} would be equivalent for local $\mathbb Q$-Gorenstein  $\mathbb C$-affine domains. If we conjecture moreover that $\mathcal B$-closure is the same as generic tight closure (as  plus closure is expected to be the same as tight closure), \eqref{i:lt}--\eqref{i:Breg} would be equivalent. Without these assumptions, it is not hard to show that if $R$ is weakly $\mathcal B$-regular, then any \qfrob\ is pure. The fact that we allow in the definition of \qreg{ity} any \qfrob, and not just powers of the non-standard Frobenius, causes an obstruction in proving that \eqref{i:uFreg}$\implies$\eqref{i:nsreg}.
\end{remark}

The importance of this new characterization of log-terminal singularities in Theorem~\ref{T:MT} is the fact that unlike the first two properties, \qreg{ity} is easily proved to descend under (cyclically) pure \homo{s}.

\begin{proposition}\label{P:pure}
Let $R\to S$ be a cyclically pure \homo\ between $\mathbb C$-affine algebras. If $S$ is \qreg, then so is $R$. 
\end{proposition}
\begin{proof}
Since $S$ is in particular normal, so is $R$ (see for instance \cite[Theorem 4.7]{SchRatSing}). Therefore, the embedding $R\to S$ is pure, by \cite[Theorem 2.6]{HoPure}. Let $c$ be a non-zero element of $R$. By assumption, we can find an \qfrob\ $\uf$ such that the $S$-module morphism 
	\begin{equation*}
	\multf c\uf S\colon S\to \bc\uf\ul S\colon x\mapsto c\uf(x)
	\end{equation*}
is pure, and whence so is its composition with $R\to S$. However, this composite morphism factors as $\multf c\uf R$ followed by the inclusion $\bc\uf\ul R \subset \bc\uf\ul S$. Therefore, the first factor, $\multf c\uf R$ is already pure, showing that $R$ is \qreg. 
\end{proof}

\begin{remark}\label{R:bou}
Proposition~\ref{P:pure} together with Theorem~\ref{T:MT} and Remark~\ref{R:bf} therefore yield the first assertion of Theorem~\ref{T:ltpure}. The last assertion   is then a direct consequence of this (after  localization), since the hypotheses imply that the inclusion $A^G\subset A$ is cyclically pure (in fact even split), where $A$ is the affine coordinate ring of $X$ and $A^G$ the subring of $G$-invariant elements (so that $X/G=\op{Spec}A^G$).
\end{remark}

\begin{remark}\label{R:reg}
Under the assumption that $S$ is regular in Theorem~\ref{T:ltpure} rather than just having log-terminal singularities, we can still conclude that $R$ has log-terminal singularities, without having to rely on the deep result by \name{Hara}. Namely, since $S$ is regular, it  is \qreg\ (Remark~\ref{R:regqreg}), whence so is $R$ by Proposition~\ref{P:pure}, and therefore $R$ has log-terminal singularities by (\ref{i:uFreg})$\implies$(\ref{i:lt}) in Theorem~\ref{T:MT}.
\end{remark}

\subsection{Log-canonical singularities}
The previous results raise the following question. \textsl{If  the non-standard Frobenius $\ulfrob\colon R\to \ul R$ is pure, for $R$ a $\mathbb Q$-Gorenstein local $\mathbb C$-affine domain, does $R$ have log-canonical singularities? Is the converse also true? What if we only require that some \qfrob\ is pure?} Note that F-pure type implies log-canonical singularities by \cite[Corollary 4.4]{Wat}, and this former condition is supposedly the analogue of $\ulfrob$ being pure. If the question and its converse  are both answered in the affirmative, we also have a positive solution to the following question: \textsl{if $R\to S$ is a cyclically pure \homo\ of $\mathbb Q$-Gorenstein local $\mathbb C$-affine domains and if $S$ has log-canonical singularities, does then so have $R$?} See also Remark~\ref{R:D0} below for some related issues.

\section{Vanishing of Maps of Tor}\label{s:vantor}

We start with providing a proof of Theorem~\ref{T:vantor} from the introduction. To this end, we need to review some results from \cite{SchBCM} on the canonical construction of big \CM\ algebras. For $R$ a local $\mathbb C$-affine domain, let $\BCM R$ be the ultraproduct of the absolute integral closures $(\seq Rp)^+$, where $\seq Rp$ is some \sr\ of $R$. We showed in \cite[Theorem A]{SchBCM} that $\BCM R$ is a (balanced) big \CM\ algebra of $R$. It follows that if $R$ is regular, then $R\to \BCM R$ is flat (\cite[Corollary 2.5]{SchBCM}). 

This construction is weakly functorial in the sense that given any local \homo\ $R\to S$ of local $\mathbb C$-affine domains, we can find a (not necessarily unique) \homo\ $\BCM R\to \BCM S$ making the following diagram commute
	\commdiagram R{} S {} {} {\BCM R} {} {\BCM S.}
If $R\to S$ is finite, then $\BCM R=\BCM S$ (see \cite[Theorem 2.4]{SchBCM}). As already observed in Remark~\ref{R:genreg}, we have the following purity result.

\begin{proposition}\label{P:cyc}
If a local $\mathbb C$-affine domain $R$ is \qreg, then $R\to \BCM R$ is cyclically pure.
\end{proposition}

Theorem~\ref{T:vantor} is a special case of the next result in view of Theorem~\ref{T:MT} and Proposition~\ref{P:pure}. It generalizes \cite[Theorem 4.12]{HHbigCM2} (I will only deal with the Tor functor here; the more general form of loc.\ cit., can be proved by the same arguments).

\begin{theorem}[Vanishing of maps of Tor]\label{T:vmt}
Let $R\to S$ be  a \homo\ of $\mathbb C$-affine algebras such that $S$ is an \qreg\ domain. Let $A$ be a regular subring of $R$ over which $R$ is module finite. Then for every $A$-module $M$ and every $i\geq 1$, the natural morphism $\tor AiMR\to \tor AiMS$ is zero.
\end{theorem}
\begin{proof}
If the map is non-zero, then it remains so after a suitable localization of $S$, so that we may assume that $S$ is local. We then may localize $A$ and $R$ at the respective contractions of the maximal ideal of $S$, and assume that $A$ and $R$ are already local. Let $\pr$ be a minimal prime of $R$ contained in the kernel of the \homo\ $R\to S$. The composition $R\to R/\pr\to S$ induces a factorization $\tor AiMR\to \tor AiM{R/\pr}\to\tor AiMS$. Thus, in order to prove the statement, it suffices to show that the second \homo\ is zero, so that we may assume that $R$ is a domain.

We have a commutative diagram
	\begin{equation}\label{d:tor}
	\CD
	A	@>>>	R	@>>>	S\\
@VVV		@VVV		@VVV	\\
	\BCM A	@= \BCM R	@>>>	\BCM S.
\endCD
	\end{equation}
Let $\phi$ be the composite morphism 
	\begin{equation}\label{eq:phi}
	\tor AiMR\to \tor AiM{\BCM R}\to \tor AiM{\BCM S}.
	\end{equation} 
By the preceding discussion and our assumptions, $\BCM A=\BCM R$ and $\BCM A$ is flat over $A$. Therefore, the middle module in \eqref{eq:phi} is zero, whence so is $\phi$. Using the commutativity of \eqref{d:tor}, we see that $\phi$ also factors as 
	\begin{equation*}
	\tor AiMR\to \tor AiMS\to \tor AiM{\BCM S}. 
	\end{equation*}
By Proposition~\ref{P:cyc}, the embedding $S\to \BCM S$ is cyclically pure, whence pure, by \cite[Theorem 2.6]{HoPure} and the fact that $S$ is   normal. Therefore, $\tor AiMS\to\tor AiM{\BCM S}$ is injective (see   \cite[Lemma 2.1(h)]{HHbigCM2}). It follows from $\phi=0$ that then also $\tor AiMR\to \tor AiMS$ must be zero, as required.
\end{proof}

The next two results follow already from \name{Hara}'s characterization of log-terminal singularities (equivalence (\ref{i:lt})$\,\Longleftrightarrow\,$(\ref{i:lt}') in Theorem~\ref{T:MT}) together with the pertinent facts on zero \ch\ tight closure. However, we can give more direct proofs using our methods. The first of these is a \BS\ type result. For regular rings, it was first proved in \cite{LT}; a tight closure proof was given in \cite{HHTC}. 

\begin{theorem}[\BS]\label{T:BS}
Let $R$ be a local   domain essentially of finite type over a field of \ch\ zero and assume $R$ has at most log-terminal singularities (or, more generally,  is \qreg). If $I$ is an ideal in $R$ generated by at most $n$ elements, then the integral closure of $I^{n+k}$ is contained in $I^{k+1}$, for every $k\geq0$.
\end{theorem}
\begin{proof}
The proof is an immediate consequence of Proposition~\ref{P:cyc} applied to \cite[Theorem B]{SchBCM}. For the reader's convenience, we repeat the argument. Let $R$ be \qreg\ and let $I$ an ideal generated by at most $n$ elements. Let $z$ be an element in the integral closure of $I^{n+k}$, for some $k\in\nat$. Take \sr{s} $\seq Rp$, $\seq Ip$ and $\seq zp$ of $R$, $I$ and $z$ respectively. Since $z$ satisfies an integral equation
	\begin{equation*}
	z^n +a_1z^{n-1}+\dots+a_n=0
	\end{equation*}
with $a_i\in I^{(n+k)i}$, we have for almost all $p$ an equation
	\begin{equation*}
	(\seq zp)^n +\seq{a_1}p(\seq zp)^{n-1}+\dots+\seq{a_n}p=0
	\end{equation*}
with $\seq{a_i}p\in(\seq Ip)^{(n+k)i}$ an \sr\ of $a_i$. In other words,  $\seq zp$ lies in the integral closure of $(\seq Ip)^{n+k}$, for almost all $p$. By \cite[Theorem 7.1]{HHbigCM2},  almost all $\seq zp$ lie in $(\seq Ip)^{k+1}\seq Rp^+\cap\seq Rp$. Taking ultraproducts, we get that $z\in I^{k+1}\BCM R\cap R$. By Proposition~\ref{P:cyc}, we get that $z\in I^{k+1}$ as required.
\end{proof}

Recall that the \emph{symbolic power} $I^{(n)}$ of an ideal $I$ in a ring $R$ is by definition the collection of all $a\in R$ for which there exists an $R/I$-regular element $s\in R$ such that $sa\in I^n$. We always have an inclusion $I^n\subset I^{(n)}$. If $I=\pr$ is prime, then $\pr^{(n)}$ is just the $\pr$-primary component of $\pr^n$. The following generalizes the main results of \cite{ELS} and \cite{HHComp} to log-terminal singularities.

\begin{theorem}\label{T:SP}
Let $R$ be a log-terminal (or, more generally,   \qreg) $\mathbb C$-affine domain. Let $\id$ be an ideal in $R$ and let $h$ be the largest height of an associated prime of $\id$ (or more generally, the largest analytic spread of $\id R_\pr$, for $\pr$ an associated prime of $R$). If $\id$ has finite projective dimension, then $\id^{(hn)}\subset \id^n$, for all $n$.
\end{theorem}
\begin{proof}
The same argument that deduces \cite[Theorem 3.4]{SchSymPow} from its positive \ch\ counterpart \cite[Theorem 1.1(c)]{HHComp}, can be used to obtain the zero \ch\ counterpart of \cite[Theorem 1.1(b)]{HHComp}, to wit, the fact that $\id^{(hn)}$ lies in the generic tight closure of $\id^n$ (use \cite[Proposition 6.3]{SchBC} in conjunction with the techniques from \cite[\S4]{SchNSTC}). By Theorem~\ref{T:MT} and Proposition~\ref{P:qreg}, each ideal is equal to its generic tight closure, proving the assertion.
\end{proof}

\section{Polarizations}\label{s:pol}

In this section, $X$ denotes a projective scheme  of finite type over some \acf\  $K$. Given an ample invertible $\loc_X$-module $\mathcal P$, we will call the pair $(X,\mathcal P)$ a \emph{polarized scheme}  and we call $\mathcal P$ a \emph{polarization} of $X$. 

Fix a polarized scheme $(X,\mathcal P)$. For each $\loc_X$-module $\mathcal F$, we define its \emph{polarization} to be the sheaf
	\begin{equation*}
	\polar{\mathcal F} := \bigoplus_{n\in\zet} \mathcal F\tensor_{\loc_X}\mathcal P^n.
	\end{equation*}
In particular, for each $\mathcal F$, we have an isomorphism
	\begin{equation*}
	\polar{\mathcal F}\iso  \mathcal F\tensor_{\loc_X}\polar{\loc_X}.
	\end{equation*}

\begin{definition}
The \emph{section ring} $S$ of $(X, \mathcal P)$ is the ring of global sections of $\polar{\loc_X}$, that is to say,
	\begin{equation}\label{eq:polsec}
	S:=\bigoplus_{n\in\zet} H^0(X,\mathcal P^n).
	\end{equation}
\end{definition}

Note that $S$ is  a finitely generated graded  algebra over $K=H^0(X,\loc_X)$ by letting $\grad Sn:=H^0(X,\mathcal P^n)$ (ampleness is used to guarantee that it is finitely generated). In fact, $S$ is positively graded, since $\mathcal P^n$ has no global sections for $n<0$. 

The polarization can be completely recovered from the section ring $S$ by the rules 
	\begin{equation*}
	X\iso \op{Proj} S \qquad\text{and}\qquad \mathcal P\iso\widetilde{S(1)}.
	\end{equation*}
In fact,  $\mathcal P^n=\widetilde{S(n)}$, for any $n\in\zet$. Global properties of $X$ can now be studied via local properties of $S$ (or more accurately, of $S_\maxim$, where $\maxim$ is the irrelevant maximal ideal generated by all homogeneous elements of positive degree).

\begin{definition}
The \emph{section module} of an $\loc_X$-module $\mathcal F$ (with respect to the polarization $\mathcal P$) is the module of global sections $H^0(X,\polar\sheaf)$ of $\polar{\mathcal F}$ and is denoted $\sect {\mathcal P}\sheaf$, or just $\sect{}\sheaf$, if the polarization is clear.
\end{definition}

In particular, the section module $\sect{}{\loc_X}$ of $\loc_X$ is just $S$ itself. Let $F:=\sect{}\sheaf$. We make $F$ into a  $\zet$-graded $S$-module by $\grad Fn:= H^0(X,\sheaf\tensor\mathcal P^n)$, for $n\in\zet$. Indeed, for each $m,n\in\zet$, we have $\grad Sm\cdot\grad Fn\subset \grad F{m+n}$ because we have canonical isomorphisms $\mathcal P^m\tensor(\sheaf\tensor\mathcal P^n)\iso \sheaf\tensor\mathcal P^{m+n}$. If $\sheaf$ is coherent, then  there is some $n_0$ such that $\sheaf\tensor\mathcal P^n$ is generated by its global sections for all $n\geq n_0$, since $\mathcal P$ is ample. Therefore, if $\sheaf$ is coherent,  $F$ is finitely generated. For each $n\in\zet$, we have an isomorphism
	\begin{equation}\label{eq:tw}
	\widetilde {\sect{}\sheaf(n)}\iso \sheaf\tensor \mathcal P^n.
	\end{equation}
Indeed, it suffices to prove that both sheaves have the same sections on each open $D_+(x)$ defined by some homogeneous element $x$, and this is straightforward. Note that we can in particular recover $\sheaf$ from its section module $\sect{}\sheaf$ since $\sheaf\iso \widetilde {\sect{}\sheaf}$.

\subsection{\Cech\ Cohomology and Polarizations}
Let $(X,\mathcal P)$ be a polarized scheme with section ring  $S:=\sect{}{\loc_X}$. Let $\maxim$ be the maximal irrelevant ideal of $S$ and let $\mathbf x=\rij xs$ be a  homogeneous system of parameters of $S$ (so that $\mathbf xS$ is in particular $\maxim$-primary). For each tuple $\mathbf i$ of indices given by $1\leq i_1<i_2<\dots<i_m\leq s$, set 
	\begin{equation*}
	\mathbf x_{\mathbf i}:= x_{i_1}x_{i_2}\cdots x_{i_m}\qquad\text{and}\qquad U_{\mathbf i}:= D_+(\mathbf x_{\mathbf i})
	\end{equation*}
(recall that $D_+(y)$ is the open consisting of all homogeneous primes of $S$ not containing the homogeneous element $y$). Let $\mathfrak U_{\mathbf x}$ be the affine open covering of $X$ given by the $U_i:=D_+(x_i)$. 

\subsubsection{\Cech\ complex of a sheaf}\label{ss:ccs}
 Let $\sheaf$ be a quasi-coherent $\loc_X$-module.  The   \emph{\Cech\ complex} of $\sheaf$ with respect to the covering $\mathfrak U_{\mathbf x}$ is the complex $\cc(\mathfrak U_{\mathbf x};\polar{\sheaf})$ given as 
	\begin{equation*}
	0\to \mathcal C^1:= \bigoplus_i H^0(U_i,\polar{\sheaf}) \to \dots\to \mathcal C^m:= \bigoplus_{\mathbf i} H^0(U_{\mathbf i},\polar{\sheaf})\to \dots
	\end{equation*}
where in $\mathcal C^m$  the index $\mathbf i$ runs over all $m$-tuples  of indices   $1\leq i_1<i_2<\dots<i_m\leq s$ and where the morphisms are, up to sign, given by restriction  (see \cite[Chapt.\ III.\ \S4]{Hart} for more details).  Using \eqref{eq:tw}, we see that $H^0( U_{\mathbf i}, \sheaf\tensor\mathcal P^n)$ is isomorphic to $\grad {\sect{\mathbf x_{\mathbf i}}\sheaf}n$. Therefore, we have a (degree preserving) isomorphism 
	\begin{equation}\label{eq:Dx}
	H^0(U_{\mathbf i},\polar{\sheaf})\iso \sect{\mathbf x_{\mathbf i}}\sheaf.
	\end{equation}

\subsubsection{\Cech\ complex of a module} 
More generally, we associate to  an arbitrary $S$-module  $F$  a \emph{\Cech\  complex} $\cc(\mathbf x;F)$ given as  
	\begin{equation}\label{eq:Ci}
	0 \to \mathcal C^1:=\bigoplus_i F_{x_i}\to \dots\to \mathcal C^m:= \bigoplus_{\mathbf i} F_{\mathbf x_{\mathbf i}}\to \dots
	\end{equation}
(with notation as above) where the morphisms are, up to sign, the natural inclusions (in fact, this construction can also be made in the non-graded case, with $\mathbf x$ an arbitrary tuple of elements in $S$; see \cite[p.\ 129]{BH} for more details). For an $\loc_X$-module $\sheaf$, we get using \eqref{eq:Dx},   an isomorphism of complexes
	\begin{equation}\label{eq:ccpolar}
	\cc(\mathfrak U_{\mathbf x};\polar{\mathcal F})\iso \cc(\mathbf x;\sect{}\sheaf).
	\end{equation}
In particular, the cohomology of either complex can be used to compute the sheaf cohomology of $\sheaf$. 

\subsubsection{Local cohomology} 
On the other hand, for arbitrary $F$, the complex $\cc(\mathbf x;F)$ can also be used to calculate local cohomology. Recall that $H^0_\maxim(F)$ is equal to the $\maxim$-torsion $\Gamma_\maxim(F)$ of $F$, that is to say, equal to the (homogeneous) submodule of all elements of $F$ which are annihilated by some power of  $\maxim$; the derived functors of $\Gamma_\maxim(\cdot)$ are then the local cohomology modules $H^i_\maxim(\cdot)$. By \cite[Theorem 3.5.6]{BH}, the local cohomology of an $S$-module $F$ can be computed as the cohomology of the augmented complex $0\to F\to \cc(\mathbf x;F)$ (that is, we inserted an additional term $\mathcal C^0:=F$ in \eqref{eq:Ci}). Hence for $i>1$,  we have an isomorphism 
	\begin{equation}\label{eq:lci}
	H^i(\cc(\mathbf x;F)) \iso H^i_\maxim(F),
	\end{equation}
whereas for $i=1$, we have a short exact sequence
	\begin{equation}\label{eq:lc1}
	0\to H^0_\maxim(F) \to F \to H^1(\cc(\mathbf x;F))\to H^1_\maxim(F)\to 0.
	\end{equation}
Since all local cohomology modules are Artinian, $F$ and $H^1(\cc(\mathbf x;F))$ have the same localizations at the various $x_i$. Moreover, using \eqref{eq:secmod} below, we get that $H^1(\cc(\mathbf x;F))=\sect{}{\tilde F}$, for $F$ a finitely generated graded $S$-module and $\tilde F$ the $\loc_X$-module associated to $F$. In conclusion, we have an equality of complexes
	\begin{equation}\label{eq:cm1}
	\cc(\mathbf x;F)=\cc(\mathbf x;\sect{}{\tilde F})
	\end{equation}
and \eqref{eq:lc1} becomes the exact sequence (see also \cite[Theorem A4.1]{Eis})
	\begin{equation}\label{eq:shf}
	0\to H^0_\maxim(F) \to F \to \sect{}{\tilde F}\to H^1_\maxim(F)\to 0.
	\end{equation}

\subsubsection{Comparison of cohomology}\label{ss:pol}
Let us summarize some of these observations. By \eqref{eq:ccoh} and \eqref{eq:ccpolar}, we have for each $i\geq0$, isomorphisms of graded $S$-modules
	\begin{equation}\label{eq:loccoh}
	H^i(X,\polar{\sheaf}) \iso H^{i+1}(\cc(\mathfrak U_{\mathbf x};\polar{\sheaf})) \iso H^{i+1}(\cc(\mathbf x; \sect{}\sheaf)). 
	\end{equation}
Moreover, for $i\geq1$, these modules are also isomorphic to $H^{i+1}_\maxim (\sect {}\sheaf)$ by \eqref{eq:lci}.  In particular, with $i=0$, isomorphism~\eqref{eq:loccoh}  becomes
	\begin{equation}\label{eq:secmod}
	\sect{}\sheaf=H^0(X,\polar{\mathcal F}) \iso   H^1(\cc(\mathfrak U_{\mathbf x};\polar{\mathcal F})).
	\end{equation}
Using that the isomorphisms in \eqref{eq:loccoh} preserve degree, we have for each $i\geq0$ and each $n\in\zet$, isomorphisms
	\begin{equation}\label{eq:loccohn}
	H^i(X, \mathcal F\tensor\mathcal P^n) \iso H^{i+1}(\grad{\cc(\mathbf x; \sect{}\sheaf)}n)\iso \grad{H^{i+1}_\maxim (\sect{}\sheaf)}n
	\end{equation}
(where the final isomorphism only holds for $i>0$).

\begin{lemma}\label{L:tens}
Let $(X,\mathcal P)$ be a polarized scheme with section ring $S$ and let $\sheaf$ and $\mathcal G$ be two  coherent $\loc_X$-modules. If $\mathcal P$ is very ample, then  there is a short exact sequence (of  degree preserving  morphisms)
	\begin{equation}\label{eq:tens}
	0\to H^0_\maxim(\sect {}\sheaf\tensor \sect{}{\mathcal G}) \to \sect {}\sheaf\tensor \sect{}{\mathcal G} \to \sect{}{\sheaf\tensor\mathcal G} \to H^1_\maxim(\sect {}\sheaf\tensor \sect{}{\mathcal G})\to 0.
	\end{equation}
\end{lemma}
\begin{proof}
Let $F:=\sect{}\sheaf$ and $G:=\sect{}{\mathcal G}$ be the respective section modules of $\sheaf$ and $\mathcal G$. By \eqref{eq:shf}, it suffices to show that the $\loc_X$-module $\widetilde{F\tensor G}$ associated to $F\tensor G$ is isomorphic with $\sheaf\tensor\mathcal G$. Since $\mathcal P$ is very ample, $S$ is generated by its linear forms and it suffices to check that both sheaves agree over each open $D_+(x)$ with $x$ homogeneous of degree one. To this end, we get, using \eqref{eq:tw},  isomorphisms
	\begin{align*}
	 (\mathcal F\tensor_{\loc_X}\mathcal G)(D_+(x)) &\iso  \mathcal F(D_+(x)) \tensor_{\loc_X(D_+(x))}  (\mathcal G)(D_+(x)) \\
	&\iso \grad {F_x}0 \tensor_{\grad {S_x}0} \grad {G_x}0\\
	&\iso \grad{(F\tensor_SG)_x}0\\
	&\iso \widetilde{(F\tensor_S G)}(D_+(x))
	\end{align*}
where  the penultimate isomorphism follows from \cite[II. Proposition 2.5.13]{EGA} (to apply this, it is necessary that $x$ has degree one). 
\end{proof}

We want to remind the reader of the following observation made in \cite{SmFano}.

\begin{proposition}\label{P:CM}
For $S$ the section ring of  a polarized scheme  $(X,\mathcal P)$, we have that $S$ is \CM\ \iff\ $H^i(X,\mathcal P^n)=0$ for all $n$ and all $0<i<\op{dim}X$. 

Under the additional assumption that  $X$ is \CM, $S$ is \CM\ \iff\ $H^i(X,\loc_X)=0$ for all $0<i<\op{dim}X$.
\end{proposition}
\begin{proof}
Let $\maxim$ be the maximal irrelevant ideal of $S$. As explained in \cite[Proposition 2.1]{HySm},  the local cohomology groups $H^0_\maxim(S)$ and $H^1_\maxim(S)$ always vanish. By a theorem of Grothendieck (\cite[Theorem 3.5.7]{BH}), $S$ is \CM\ \iff\ $H^{i+1}_\maxim(S)=0$, for all $i<\op{dim} X$. By \eqref{eq:loccoh} in \S\ref{ss:pol}, this in turn is equivalent with $H^i(X,\mathcal P^n)=0$, for all $0<i<\op{dim}X$ and all $n$, proving the first assertion.

Suppose $X$ is moreover \CM. Since $\mathcal P$ is invertible, $H^i(X,\mathcal P^n)=0$ for all $0<i<\op{dim}X$ and $n\ll 0$ by Serre duality (\cite[III. Theorem 7.6]{Hart}). The same is true for $n\gg 0$, since $\mathcal P$ is ample (\cite[III. Proposition 5.3]{Hart}). Therefore polarizing $X$ with respect to a sufficiently large power $\mathcal P^s$ instead of $\mathcal P$, we may even assume that $H^i(X,\mathcal P^n)=0$ for all $n\neq 0$ and $0<i<\op{dim}X$. The second assertion then follows from this and the first assertion (applied to the new polarization).
\end{proof}

\subsection{Absolute Frobenius}
Assume now that $X$ is a scheme of finite type over a perfect field of prime \ch\ $p$. Fix a power $q$ of $p$. We call the endomorphism  which is the identity on the underlying topological space and the $q$-th power map $x\mapsto x^q$ on the structure sheaf, an (absolute) \emph{Frobenius} on $X$ and we also denote it $\frob q$. For instance, if $X=\op{Spec} C$, then $\frob q\colon X\to X$ is the endomorphism of $X$ determined by the (ring-theoretic) Frobenius $\frob q\colon C\to C$. Recall that to emphasize the $C$-algebra structure on $C$ induced by $\frob q$, we write $\bc{\frob q}C$. In particular, the natural map  $C\to \bc{\frob q}C$ is just $\frob q$ and induces the absolute Frobenius on $X=\op{Spec}C$.

If $\mathcal L$ is an invertible sheaf on $X$, then $\frob q^*\mathcal L\iso \mathcal L^q$. Indeed, this is true locally on an affine open cover, and hence also globally. In particular, if $\mathcal P$ is a polarization on $X$, then so is $\frob q^*\mathcal P\iso\mathcal P^q$. If $S$ is the section ring of $(X,\mathcal P)$, then $S^{(q)}$ is the section ring of $(X,\mathcal P^q)$, where in general for a graded module $M$, we write $M^{(a)}$ to denote the graded module $\oplus_n \grad M{an}$. Since $\frob q$ on $S$ is just the ring \homo\ $S\to S^{(q)}$ given by $x\mapsto x^q$, we see that $\bc{\frob q}S\iso S^{(q)}$. Taking projective spaces, we get the absolute Frobenius $X=\op{Proj}(S^{(q)})\to X=\op{Proj}(S)$ (here we used \eqref{eq:polsec} twice).

\section{Vanishing Theorems}

In \cite{SmFano}, Smith introduces the notion of a globally $F$-regular type variety and shows that it admits several vanishing theorems. She moreover conjectures that any GIT quotient by a reductive group of a smooth Fano variety satisfies these vanishing theorems.  We will establish this conjecture using the notion of \qreg{ity} as a substitute for F-regularity.

Let $X$ be a connected projective scheme of finite type over $\mathbb C$ (a \emph{projective variety}, for short).

\begin{definition}
We say that $X$ is \emph{globally \qreg}, if some section ring $S$ of $X$ is \qreg\ at its vertex, that is to say, $S_\maxim$ is \qreg\ where $\maxim$ is the  irrelevant maximal ideal of $S$.
\end{definition}

\begin{remark}
In  \cite{SmFano}, \name{Smith} calls $X$ \emph{globally F-regular}  if some section ring $S$   is strongly F-regular type (note that since $S$ is positively graded, this is equivalent  by \cite{LyuSm} with $S$ being weakly  F-regular type). By Proposition~\ref{P:qreg}, this implies that $S$ is \qreg\ and hence that $X$ is globally \qreg. Moreover, if $S_\maxim$ is ($\mathbb Q$-)Gorenstein, where $\maxim$ is the maximal irrelevant ideal, then globally F-regular type and globally \qreg\ are equivalent in view of  Theorem~\ref{T:MT}.

So we could deduce the desired vanishing theorems from the work of \name{Smith} in \cite{SmFano}, if we are willing to use \name{Hara}'s characterization of F-regular type. However, using a non-standard version of her arguments, we can as easily derive the vanishing theorems directly, without  any appeal to \name{Hara}'s work (and hence without using Kodaira Vanishing).
\end{remark} 

\begin{remark}\label{R:all}
As in \cite{SmFano}, one can prove directly that if $X$ is globally \qreg, then every section ring is locally \qreg\ at its vertex. Alternatively, this follows from \cite[Theorem 3.10]{SmFano} (even without localizing at the irrelevant maximal ideal), if we use Theorem~\ref{T:MT} as in the previous remark. 

In that respect, note that if the section ring $S$ with respect to the polarization $\mathcal P$ is \qreg\ at its vertex, then so is any Veronese subring $S^{(r)}$ by Proposition~\ref{P:pure}, as it is a pure subring. In particular, any section ring corresponding to a positive power of $\mathcal P$ is \qreg\ at its vertex. In particular, we may always assume, without relying on the results from \cite{SmFano},  that a globally \qreg\ variety admits a \emph{very} ample polarization whose section ring is \qreg\ at its vertex.
\end{remark}

As already mentioned, the main advantage of using \qreg{ity} instead of F-regular type is the fact that it descends under pure \homo\ (Proposition~\ref{P:pure}). In particular, we get the following descent property for quotient singularities.

\begin{theorem}\label{T:GIT}
Let $X$ be a connected projective variety over $\mathbb C$. Let $G$ be a reductive group acting algebraically on $X$ and let $X/\!/G$ be an arbitrary GIT quotient of $X$. If $X$ is globally \qreg, then so is $X/\!/G$.
\end{theorem}
\begin{proof}
Any GIT quotient of $X$ is obtained by taking some polarization $\mathcal P$ of $X$, extending the $G$-action to $\mathcal P$, taking the section ring $S$ with the induced $G$-action and letting $X/\!/G:=\op{Proj}S^G$, where $S^G$ is the ring of invariants of $S$. In particular, $S^G$ is a section ring of $X/\!/G$. Since $S$ is \qreg\ by Remark~\ref{R:all}, so is $S^G$ by Proposition~\ref{P:pure} as $S^G\subset S$ is pure (even split).
\end{proof}

\begin{theorem}\label{T:van}
Let $X$ be a globally \qreg\ connected projective variety over $\mathbb C$ and let $\mathcal F$ be an invertible $\loc_X$-module. If for some $i>0$ and some effective Cartier divisor $D$, all $H^i(X,\mathcal F^n(D))$ vanish  for $n\gg 0$, then $H^i(X,\mathcal F)$ vanishes.
\end{theorem}
\begin{proof}
Choose a polarization $\mathcal P$ of $X$ with section ring $S$, so that $S_\maxim$ is \qreg, where $\maxim$ is the irrelevant maximal ideal of $S$. By Remark~\ref{R:all}, we may assume without loss of generality that $\mathcal P$ is very ample. Let  $I$ be the section module of $\mathcal I:=\loc_X(D)$. Let $\mathbf x$ be a homogeneous system of parameters of $S$ and let $\mathfrak U_{\mathbf x}$ be the open affine covering given by the $D_+(x_i)$. Since $D$ is Cartier, $I$ is a fractional ideal, that is to say, a finitely generated rank-one $S$-submodule of the field of fractions $K$ of $S$. Clearing denominators in the inclusion $I\subset K$, we can find an $S$-module morphism $\psi\colon I\to S$. Since $D$ is effective, $I$ admits a canonical section   $s\in\grad I0=H^0(X,\mathcal I)$ (see for instance \cite[B.4.5]{Ful}). In particular, the morphism $S\to I\colon 1\mapsto s$ is degree preserving. Put $c:=\psi(s)$. By \qreg{ity}, there is an \qfrob\ $\uf$ such that $\multf c\uf{S_\maxim}$ is pure. The composition $S\to I\map\psi S\colon 1\mapsto s\mapsto c$  is equal to multiplication with $c$ on $S$ (where we disregard the grading). Tensoring this composite \homo\  with $\bc\uf\ul S$ gives 
	\begin{equation*}
	 \bc\uf\ul S\to I\tensor\bc\uf\ul S\to \bc\uf\ul S\colon 1\mapsto s\tensor 1\mapsto c
	\end{equation*} 
which composed with the inclusion $S\into \bc\uf\ul S$ therefore gives the morphism $\multf c\uf S$. By assumption, the base change $\multf c\uf {S_\maxim}$  is pure. Since
	\begin{equation*}
	S_\maxim\to S_\maxim\tensor I\tensor \bc\uf\ul S
	\end{equation*}
is a factor of the pure morphism $\multf c\uf {S_\maxim}$, it is also pure. Let $F:=\sect{}\sheaf$ be the section module of $\mathcal F$. Tensoring with $F$ yields a pure $S_\maxim$-module morphism
	\begin{equation*}
	F_\maxim\to  (F\tensor I\tensor \bc\uf\ul S)_\maxim.
	\end{equation*}
Using the isomorphism  $\bc\uf F\iso F\tensor\bc\uf\ul S$ (see \S\ref{s:uf}), we can identify $F\tensor I\tensor\bc\uf\ul S$ with $I\tensor\bc\uf F$. Taking \Cech\ complexes with respect to the tuple $\mathbf x$ yields a pure homomorphism of \Cech\ complexes
	\begin{equation*}
	\cc(\mathbf x;F_\maxim)   \to \cc(\mathbf x;  (I\tensor\bc\uf F)_\maxim). 
	\end{equation*}
It is well-known that purity is preserved after taking cohomology, so that we have a pure morphism
	\begin{equation*}
	H^{i+1}(\cc(\mathbf x;F_\maxim)) \into H^{i+1} (\cc(\mathbf x;  (I\tensor\bc\uf F)_\maxim)).
	\end{equation*}
Since at a maximal ideal, the cohomology of a localized \Cech\ complex is the same as the cohomology of the non-localized \Cech\ complex (see for instance \cite[Remark 3.6.18]{BH}), we get an injective morphism
	\begin{equation}\label{eq:HL}
	H^{i+1}(\cc(\mathbf x;F))\into H^{i+1} (\cc(\mathbf x;  I\tensor\bc\uf F)) 
	\end{equation}
By a similar argument as in \S\ref{s:grad}, each module in $\cc(\mathbf x;  I\tensor\bc\uf F)$, although not  graded, has a graded piece in each (non-standard) degree. This property is inherited by the cohomology groups and \eqref{eq:HL} preserves degrees. Hence  in degree zero, we get an injective morphism 
	\begin{equation}\label{eq:HL0}
	\grad{H^{i+1}(\cc(\mathbf x;F))}0\into \grad{H^{i+1} (\cc(\mathbf x;  I\tensor\bc\uf F))}0.
	\end{equation} 
I claim that the   right hand side of \eqref{eq:HL0} is zero, whence by injectivity, so is the  left hand side. Since the latter is just $H^i(X,\mathcal F)$ by \eqref{eq:loccohn}, the theorem follows from the claim.

To prove the claim, let $(\seq Xp,\seq{\mathcal P}p)$, $\seq Sp$, $\seq{\mathbf x}p$, $\seq{\mathcal F}p$ and $\seq{\mathcal I}p$ be \sr{s} of $(X,\mathcal P)$, $S$, $\mathbf x$, $\mathcal F$ and $\mathcal I$ respectively, and suppose the \qfrob\ $\uf$ is given as the ultraproduct of the Frobenii $\frob q$ (for $q:=p^{\seq ep}$ some power of $p$). By \S\ref{s:polar}, almost all $(\seq Xp,\seq{\mathcal P}p)$ are polarized. Using Theorem~\ref{T:nscoh},   the section modules $\seq Fp:=\sect{}{\seq\sheaf p}$ and $\seq Ip:=\sect{}{\seq{\mathcal I}p}$ are \sr{s} of respectively $F$ and $I$. In particular, the ultraproduct of the $\bc{\frob q}\seq Fp\iso\seq Fp\tensor \bc{\frob q}\seq Sp$ is equal to $\bc\uf F$ and we have an isomorphism of  \Cech\ complexes
	\begin{equation*}
	\cc(\mathbf x;I\tensor\bc\uf F)\iso \up p\cc(\seq{\mathbf x}p; \seq Ip\tensor\bc{\frob q}\seq Fp).
	\end{equation*}
Since cohomology commutes with ultraproducts, we get an isomorphism
	\begin{equation*}
	H^{i+1}(\cc(\mathbf x;I\tensor\bc\uf F))\iso \up p H^{i+1}(\cc(\seq{\mathbf x}p; \seq Ip\tensor\bc{\frob q}\seq Fp)).
	\end{equation*}
Therefore, the claim follows if we can show that almost all
	\begin{equation}\label{eq:ccup}
	\grad {H^{i+1}(\cc(\seq{\mathbf x}p; \seq Ip\tensor\bc{\frob q}\seq Fp))}0=0.
	\end{equation}

Let $\mathfrak U_{\seq{\mathbf x}p}$ be the affine covering of $\seq Xp$ given by the $D_+(\seq{x_i}p)$. Since they hold  locally on the covering  $\mathfrak U_{\seq{\mathbf x}p}$, we have isomorphisms of $\loc_{\seq Xp}$-modules
	\begin{equation*}
	\widetilde{(\bc{\frob q}\seq Fp)} \iso \bc{\frob q}\seq\sheaf p \iso \seq\sheaf p^q.
	\end{equation*}
Applying Lemma~\ref{L:tens} twice shows that $\sect{}{\seq{\mathcal I}p\tensor\seq{\mathcal F}p^q}$ and $\seq Ip\tensor\bc{\frob q}\seq Fp$ are isomorphic up to $\maxim$-torsion. In particular, this yields an isomorphism of   \Cech\ complexes
	\begin{equation*}
	\cc(\seq{\mathbf x}p; \seq Ip\tensor\bc{\frob q}\seq Fp) \iso \cc(\mathfrak U_{\seq{\mathbf x}p}; \polar{(\seq{\mathcal I}p\tensor\seq{\mathcal F}p^q)}).
	\end{equation*}
Taking  cohomology, we get 
	\begin{equation*}
	 H^{i+1}(\cc(\seq{\mathbf x}p; \seq Ip\tensor\bc{\frob q}\seq Fp))\iso H^{i+1}(\cc(\mathfrak U_{\seq{\mathbf x}p}; \polar{(\seq{\mathcal I}p\tensor\seq{\mathcal F}p^q)})).
	\end{equation*}
By \eqref{eq:loccohn} in \S\ref{ss:pol}, the zero-th homogeneous part of the right hand side is isomorphic to $H^i(\seq Xp,\seq{\mathcal I}p\tensor\seq{\mathcal F}p^q)$, and this is zero for large enough $p$  by Corollary~\ref{C:van} and our assumption. Thus, we showed \eqref{eq:ccup} and hence finished the proof.
\end{proof}

As in \cite{SmFano}, we immediately obtain the following corollaries (together with Remark~\ref{R:bf}, they prove Theorem~\ref{T:KVvan}).

\begin{corollary}\label{C:nef}
Let $X$ be a globally \qreg\ connected projective variety over $\mathbb C$ and let $\mathcal L$ be an invertible $\loc_X$-module. If $\mathcal L$ is numerically effective (NEF), then $H^i(X,\mathcal L)$ vanishes, for all $i>0$.
\end{corollary}
\begin{proof}
Suppose first that $\mathcal L$ is ample. By Serre Vanishing, $H^i(X,\mathcal L^n)=0$ for $n\gg 0$ and $i>0$. Hence $H^i(X,\mathcal L)=0$ by Theorem~\ref{T:van}. Suppose now that $\mathcal L$ is merely NEF. This means that $\mathcal L^n(D)$ is ample, for all $n\geq 0$, where $D$ is some ample effective Cartier divisor. Since we already proved the ample case, $H^i(X,\mathcal L^n(D))=0$, for all $n\geq 0$ and $i>0$. Therefore, $H^i(X,\mathcal L)=0$ by another application of Theorem~\ref{T:van}.
\end{proof}

\begin{corollary}[Kawamata-Viehweg Vanishing]\label{C:bignef}
Let $X$ be a connected projective variety over $\mathbb C$ and let $\mathcal L$ be an invertible $\loc_X$-module. If  $X$ is globally \qreg\ and if $\mathcal L$ is big and NEF, then $H^i(X,\mathcal L^{-1})=0$, for all $i<\op{dim}X$.
\end{corollary}
\begin{proof}
Fix some $i<\op{dim}X$. Because $\mathcal L$ is big and NEF, we can find an effective Cartier divisor $D$ such that $\mathcal L^m(-D)$ is ample for all $m\gg 0$, by \cite[Proposition 2.61]{KM}. Let $S$ be a section ring of $X$ which is \qreg, whence in particular \CM. Hence $X$ is in particular \CM. It follows from Proposition~\ref{P:CM} that every section ring of $X$ is then \CM\ (since the criteria only depends on $X$). In other words, as in the proof of that proposition, we get for any ample invertible sheaf $\mathcal P$ that $H^i(X,\mathcal P^n)=0$, for all $n$ of sufficiently large absolute value. Applied with $\mathcal P:=\mathcal L^m(-D)$, we get that 
	\begin{equation*}
	H^i(X,(\mathcal L^m(-D))^{-n})=H^i(X,(\mathcal L^{-m}(D))^n)=0
	\end{equation*} 
for all sufficiently large $m$ and $n$. Hence, for fixed $m$,  Theorem~\ref{T:van} yields the vanishing of $H^i(X,\mathcal L^{-m}(D))$. Since this holds for all large $m$, another application of Theorem~\ref{T:van} finally gives $H^i(X,\mathcal L^{-1})=0$.
\end{proof}

\begin{remark}\label{R:D0}
Call a $\mathbb C$-affine domain $R$  \emph\qpure, if $R\to \bc\uf{\ul R}$ is pure for some \qfrob\ $\uf$. Call a connected projective variety  $X$ over $\mathbb C$ \emph{globally \qpure}, if some section ring of $X$ is \qpure. Inspecting the proof of Theorem~\ref{T:van}, we get the following weaker version ($D=0$): \textsl{if $X$ is globally \qpure\ and $\mathcal L$ invertible with $H^i(X,\mathcal L^n)=0$ for all $n\gg0$, then $H^i(X,\mathcal L)=0$.} In particular, the argument in the proof of Corollary~\ref{C:nef} shows that on a globally \qpure\ variety, an ample invertible sheaf has no higher cohomology. 

In fact, we can even prove Kodaira Vanishing for this class of varieties: \textsl{if $X$ is globally \qpure\ and \CM, then $H^i(X,\inv{\mathcal L})=0$ for all $i<\op{dim}X$ and all ample invertible sheaves $\mathcal L$ on $X$.} Indeed, by Serre duality (\cite[III. Corollary 7.7]{Hart}), the dual of $H^i(X,\mathcal L^{-n})$ is $H^{d-i}(X, \omega_X\tensor\mathcal L^n)$ where $d$ is the dimension of $X$ and $\omega_X$ the dualizing sheaf on $X$. Because $\mathcal L$ is ample, the latter cohomology group vanishes for large $n$, and hence so does the first. Applying the weaker version of the vanishing theorem to this, we get that $H^i(X,\inv{\mathcal L})$ vanishes.

Because of the analogy with the notion of \emph{Frobenius split}  (see \cite[Proposition 3.1]{SmFano}) and the fact that a Schubert variety has this property (\cite[Theorem 2]{MR}), it is reasonable to expect that a Schubert variety is globally \qpure. As a corollary, we would obtain Kodaira Vanishing for GIT quotients of Schubert varieties, since \qpurity\ descends under pure \homo{s} (by the same argument as for Proposition~\ref{P:pure}).
\end{remark}

\section{Fano Varieties}\label{s:Fano}

Let $X$ be a connected, normal projective variety  over $\mathbb C$. The canonical (or, \emph{dualizing}) sheaf $\omega_X$ of $X$ is the unique reflexive sheaf which agrees with the sheaf $\wedge^d\Omega_{X/\mathbb C}$  on the smooth locus of $X$. We call $X$ \emph{Fano}, if its anti-canonical sheaf $\omega_X^{-1}$ is ample (we do not require $X$ to be smooth).

\begin{theorem}\label{T:Fano}
A Fano variety   over $\mathbb C$ with rational singularities is globally \qreg.
\end{theorem}
\begin{proof}
Let $X$ be a Fano variety   with rational singularities. Let $S$ be the anti-canonical section ring of  $X$, that is to say, the section ring with respect to the polarization given by the ample sheaf $\omega_X^{-1}$. It is well-known (see for instance \cite[Proposition 6.2]{SmFano}), that $S$ is Gorenstein and has again rational singularities. Since rational Gorenstein singularities are  log-terminal, we obtain from Theorem~\ref{T:MT} that $S_\maxim$ is \qreg, showing that $X$ is  globally \qreg.
\end{proof}

\begin{remark}
In proving that a Fano variety with rational singularities is globally \qreg, we have used Kodaira Vanishing twice:  via \name{Hara}'s result in Theorem~\ref{T:MT} and via \cite[Proposition 6.2]{SmFano}. Combining Theorems~\ref{T:KVvan} and \ref{T:GIT} with the previous theorem yields Theorem~\ref{T:Fanoquot} from the introduction. 
\end{remark}


\providecommand{\bysame}{\leavevmode\hbox to3em{\hrulefill}\thinspace}

\end{document}